\documentclass[12pt]{article}

\usepackage{amssymb}
\usepackage{amsmath}
\usepackage{color}
\usepackage{graphics}
\usepackage{epsfig}

\newtheorem{claim}{\bf \t}[part]

\newtheorem{Lemma}{Lemma}[part]

\newtheorem{Remark}{Remark}[part]
\newtheorem{Theorem}{Theorem}[part]

\numberwithin{Assumption}{section} \numberwithin{Corollary}{section}
\numberwithin{Definition}{section} \numberwithin{equation}{section}
\numberwithin{Example}{section} \numberwithin{Lemma}{section}
\numberwithin{Proposition}{section} \numberwithin{Remark}{section}
\numberwithin{Theorem}{section}

\def\v{\varepsilon}

\def\x{\xi}
\def\t{\theta}
\def\T{\Theta}

\def\a{\alpha}
\def\b{\beta}
\def\g{\gamma}
\def\d{\delta}
\def\l{\lambda}
\def\f{\frac}

\def\r{\rho}

\def\s{\sigma}

\def\o{\omega}
\def\di{\displaystyle}
\def\i{\infty}

\def\text#1{{\rm #1}}
\begin{document}
\date{}
\title{\Large \bf Global well-posedness of 2D compressible Navier-Stokes equations with
large data and vacuum}
\author{\small \textbf{Quansen Jiu},$^{1,3}$\thanks{The research is partially
supported by National Natural Sciences Foundation of China (No.
11171229) and Project of Beijing Education Committee. E-mail:
jiuqs@mail.cnu.edu.cn}\quad
  \textbf{Yi Wang}$^{2,3}$\thanks{The research is partially supported by Hong Kong RGC Earmarked Research Grant CUHK4042/08P and
CUHK4041/11P, and National Natural Sciences Foundation of China (No.
11171326). E-mail: wangyi@amss.ac.cn.}\quad and \textbf{Zhouping
Xin}$^{3}$\thanks{The research is partially supported by Zheng Ge Ru
Funds, Hong Kong RGC Earmarked Research Grant CUHK4042/08P and
CUHK4041/11P, and a Focus Area Grant at The Chinese University of
Hong Kong. Email: zpxin@ims.cuhk.edu.hk}} \maketitle \small $^1$
School of Mathematical Sciences, Capital Normal University, Beijing
100048, P. R. China

\small $^2$ Institute of Applied Mathematics, AMSS, and Hua Loo-Keng
Key Laboratory of Mathematics, CAS, Beijing 100190, P. R. China

\small $^3$The Institute of Mathematical Sciences, Chinese
University of HongKong, HongKong\\

 {\bf Abstract:} In this paper, we study the global well-posedness of the 2D compressible Navier-Stokes equations with
large initial data and vacuum. It is proved that if the shear viscosity $\mu$ is a positive constant and the bulk viscosity $\l$ is the power function of the density, that is, $\l(\r)=\r^\b$ with $\b>3$, then the 2D compressible Navier-Stokes equations with the periodic boundary conditions on the torus $\mathbb{T}^2$ admit a unique global classical solution $(\r,u)$ which may contain vacuums in an open set of $\mathbb{T}^2$.  Note that the initial data can be arbitrarily large to contain vacuum states.



\section{Introduction } \setcounter{equation}{0}
\setcounter{Assumption}{0} \setcounter{Theorem}{0}
\setcounter{Proposition}{0} \setcounter{Corollary}{0}
\setcounter{Lemma}{0} In this paper, we consider the following
compressible and isentropic Navier-Stokes equations with
density-dependent viscosities
\begin{eqnarray}\label{CNS}
\left\{ \begin{array}{ll}
\partial_t\rho+{\rm div}(\rho u)=0,  \\
\partial_t(\rho u) + {\rm div} (\rho u\otimes u) + \nabla P(\r)=\mu\Delta u+\nabla((\mu+\l(\r)){\rm div}u), &\quad x\in\mathbb{T}^2, t>0,
\end{array}
\right.
\end{eqnarray}
where $\rho (t, x)\geq 0 $, $u(t, x)=(u_1,u_2)(t,x) $ represent the
density and the velocity of the fluid, respectively. And
$\mathbb{T}^2$ is the 2-dimensional torus $[0,1]\times[0,1]$ and
$t\in[0,T]$ for any fixed $T>0$. We denote the right hand side of
$\eqref{CNS}_2$  by
$$
\mathcal{L}_\r u=\mu\Delta u+\nabla((\mu+\l(\r)){\rm div}u).
$$
Here, it is assumed that
\begin{equation}\label{1.2}
\mu={\rm const.}>0,\qquad \l(\r)=\r^\beta, \quad \beta>3,
\end{equation}
such that the operator $\mathcal{L}_\r$ is strictly elliptic.

 Let the pressure function be given by
\begin{eqnarray}
P(\r)=A \r^{\gamma },
\end{eqnarray}
where $\gamma >1$ denotes the adiabatic exponent and $A>0$ is the
constant. Without loss of generality, $A$ is normalized to be $1$.
The initial values are given by
\begin{equation}\label{initial-v}
(\r,u)(t=0,x)=(\r_0,u_0)(x).
\end{equation}
Here the periodic boundary conditions on the unit torus
$\mathbb{T}^2$ on $(\r,u)(t,x)$ are imposed to the system \eqref{CNS}.
This model problem, \eqref{CNS}-\eqref{initial-v}, was first proposed by Vaigant-Kazhikhov in \cite{Kazhikhov}  where they showed the well-posedness of the classical solution to this problem provided the initial density is uniformly away from vacuum. In this paper, we study the global well-posedness of the classical solution to this problem \eqref{CNS}-\eqref{initial-v} with general nonnegative initial densities.

There are extensive studies on global well-posedness of the
compressible Navier-Stokes equations in the case that both the shear
and the bulk viscosity are positive constants satisfying the
physical restrictions. In particular, the one-dimensional theory is
rather satisfactory, see \cite{Hoff-Smoller, lius, Ka, KS} and the
references therein. In  multi-dimensional case,
 the local well-posedness theory of classical solutions to both initial-value and initial-boundary-value problems was
 established  by Nash \cite{Nash}, Itaya \cite{Itaya} and Tani \cite{Tani} in the absence of vacuum. The short time well-posedness
 of either strong or classical solutions containing vacuum was studied recently by Cho-Kim \cite{CK1} and Luo\cite{Luo} in 3D and 2D case, respectively. In particular, Cho-Kim \cite{CK1}
 obtained the short existence and uniqueness of the classical solution to the Cauchy problem for the isentropic CNS
 with general nonnegative initial density under the assumption that the initial data satisfies a
 natural compatibility condition \cite{CK1}. One of the fundamental questions is whether these local (in time) solutions can be extended globally in time. The first pioneering work along this line is the well-known theory of Matsumura-Nishida \cite{MN}, where they obtained a unique global classical solution to the CNS in $H^s(\mathbb{R}^3)$ $(s\geq 3)$ for initial data close to its far field state which is a non-vacuum equilibrium state, and furthermore, the solution behaves diffusively toward the far field state. The proof in \cite{MN} consists of elaborate energy estimates based on the dissipative structure of the CNS and spectrum analysis for the linearized of CNS at the non-vacuum far field state. This theory has been generalized to data with discontinuities by Hoff \cite{Hoff1} and data in Besov spaces by Danchin in \cite{Dachin}. It should be noted that this theory
\cite{MN, Hoff1, Dachin} requires that the solution has small
oscillations from the uniform non-vacuum far field state so that the
density is strictly away from the vacuum uniformly in time. A
natural and important long standing open problem is whether a
similar theory holds for the initial data containing vacuums. In
this direction, the major breakthrough is due to P. L. Lions
\cite{Lions}, where he obtained the existence of a renormalized weak
solution with finite energy and large initial data which can contain
vacuums for the isentropic CNS when the exponent $\gamma$ is
suitably large, see also the refinements and
 generalizations in \cite{F1, JZ}. However, little is known on the structure, regularity, and
 uniqueness of such a weak solution except the partial regularity estimates for 2-dimensional periodic problems in Desjardins \cite{Des}
 where a stronger estimate is obtained under the assumption of uniform boundedness of the density. Recently, under some
 additional assumptions on the viscosity coefficients, and the far fields state is a non-vacuum
  state, Hoff \cite{Hoff1, Hoff2} obtained a new type of global weak solution with small total energy for the isentropic CNS, which
  have extra structure and regularity information (such as Lagrangian structure in the non-vacuum region) compared with
  the renormalized weak solutions in \cite{Lions, F1, JZ}. However, the uniqueness and regularity of those weak solutions whose
  existence has been proved in \cite{Lions, F1, JZ} remain completely open in general. By the weak-strong uniqueness of
  P. L. Lions \cite{Lions}, this is equivalent to the problem of global (in time) well-posedness of classical solution in the
  presence of vacuum. It should be pointed out that this important question is a very difficult and subtle issue since,
  in general, one would not expect a positive answer to this question due to the finite time blow-up results of Xin
  in \cite{Xin}, where it is shown that in the case that the initial density has compact support, any smooth solution to the
  Cauchy problem of the CNS without heat conduction blows up in finite time for any space dimension, see also the recent generalizations to the case for non-compact but rapidly decreasing (at far fields) initial density \cite{R}. The mechanism
  for such a blow-up has also been investigated recently and various blow-up criterion have been derived in \cite{FJ, FZZ, hlx1, hlx, hx2, SWZ1, SWZ}. More recently, Huang-Li-Xin\cite{hlx4} proved the global well-posedness of classical solutions with small energy but large oscillations which can contain vacuums to 3D isentropic compressible Navier-Stokes equations. See also the recent generalization to 3D full compressible Navier-Stokes equations \cite{Huang-Li}, the isentropic Navier-Stokes equations with potential forces \cite{LZZ}, and 1D or spherically symmetric isentropic Navier-Stokes equations with large initial data \cite{DWZ, DWZ1}.

The case that the viscosity coefficients depend on the density and
vanish at the vacuum has received a lot attention recently, see \cite{BD1, BDL, BD3, BDG, Dachin, GJX, J, JXZ, JZ, JWX, JWX1, JX, LLX, LXY, MV, MV2, SZ, YYZ, YZ2, ZF} and the references therein.
Liu, Xin and Yang first proposed in \cite{LXY} some
models of the compressible Navier-Stokes equations with
density-dependent viscosities to investigate the dynamics of the
vacuum. On the other hand, when deriving by Chapman-Enskog
expansions from the Boltzmann equation, the viscosity of the
compressible Navier-Stokes equations  depends on the temperature and
thus on the density for isentropic flows. Also, the viscous
Saint-Venant system for the shallow water, derived from the
incompressible Navier-Stokes equations with a moving free surface,
is expressed exactly as in \eqref{CNS} $N=2$, $\mu=\rho$, $\l=0$, and $\gamma=2$ (see \cite{GP}).
For the special case, \eqref{1.2}, the global well-posedness result of Vaigant-Kazhikhov \cite{Kazhikhov} is the first important surprising result for general large initial data with the only constraint that it is initially away from vacuum. However, in
the presence of vacuum, there
 appear new mathematical challenges in dealing with such systems. In particular, these systems
become highly degenerate. The velocity  cannot even be defined in
the presence of vacuum and hence it is difficult to get  uniform
estimates for the velocity near vacuum. Substantial achievements have been made for
the one-dimensional case, such as both short time and long time
existence and uniqueness for the problem of a compact of viscous
fluid expands into vacuum with either stress free condition or
continuity condition have been established with
$\mu=\rho^\alpha$ for suitable $\alpha$, see \cite{LXY, JX, YYZ, YZ2} etc.
 Li-Li-Xin \cite{LLX} recently proved the global existence of
weak solutions to the initial-boundary value problem for such a
system on a finite internal with general initial data which may
contain vacuum and discovered the phenomena that all the
vacuum states must vanish in finite time and any smooth solution
blows up near the time of vacuum vanishing which are in sharp
contrast to the case of constant viscosity coefficients, which have
been extended to the Cauchy problem on $\mathbb{R}^1$ for arbitrary
initial data with a uniform non-vacuum state at far fields by
Jiu-Xin \cite{JX}. In the case that a basic nonlinear wave pattern
is the rarefaction wave, whose nonlinear asymptotic stability has
been proved in \cite{JWX, JWX1} for the one-dimensional isentropic
CNS system with density-dependent viscosity in the framework of weak
solutions even the rarefaction wave is connecting to the
vacuum\cite{lius}.
Note also that in the case that the initial data
is strictly away from vacuum, Mellet and Vasseur has obtained the
existence and uniqueness of global strong solution to the
one-dimensional Cauchy problem \cite{MV2}. However, the progress is
very limited for multi-dimensional problems. Even the short time
well-posedness of classical solutions has not been established for
such a system in the presence of vacuum. The global existence of
general weak solutions to the compressible Navier-Stokes equations
with density-dependent viscosities or the viscous Saint-Venant
system for the shallow water model in the multi-dimensional case
remains open, and one can refer to \cite{BD3}, \cite{BDG},
\cite{GJX}, \cite{MV} for recent developments along this line.  Note
also that Zhang-Fang \cite{ZF} proved the existence of global weak
solution with small energy to 2D Vaigant-Kazhikhov model
\cite{Kazhikhov} in the framework of \cite{Hoff2} and presented the
vanishing vacuum behavior. However, the uniqueness of this weak
solution is open.

In this paper, we investigate the global existence of the classical solution to 2-dimensional Vaigant-Kazhikhov model \cite{Kazhikhov}, that is, CNS system \eqref{CNS}-\eqref{initial-v} with periodic boundary condition and general nonnegative
initial density. It should be noted that for the 2-dimensional problem, the basic reformulation of Vaigant-Kazhikhov \cite{Kazhikhov} and the formulation in terms of the material derivative used in \cite{Hoff1, hlx4} are equivalent. Following some of the key ideas developed by Vaigant-Kazhikhov \cite{Kazhikhov}, we are able to derive the uniform upper bound of the density under the assumptions that the initial density is nonnegative. Then we can derive the higher order estimates to the solution to guarantee the existence of the global classical solution.

The main results of the present paper can be stated in the following.

\begin{Theorem}\label{theorem2}
If the initial values $(\r_0,u_0)(x)$ satisfy that
\begin{equation}\label{in-d1}
0\leq(\r_0(x), P(\r_0)(x))\in
W^{2,q}(\mathbb{T}^2)\times W^{2,q}(\mathbb{T}^2),\quad u_0(x)\in H^2(\mathbb{T}^2),\quad \int_{\mathbb{T}^2}\r_0(x) dx>0
\end{equation}
for some $q>2$ and the compatibility condition
\begin{equation}\label{cc}
\mathcal{L}_{\r_0}u_0-\nabla P(\r_0)=\sqrt\r_0  g(x)
\end{equation}
with some $g\in L^2(\mathbb{T}^2)$, then there exists a unique
global classical solution $(\r,u)(t,x)$ to the compressible
Navier-Stokes equations \eqref{CNS}-\eqref{initial-v} with
\begin{equation}\label{Jan-16-1}
\begin{array}{ll}
\di 0\leq \r(t,x)\leq C,\quad \forall(t,x)\in [0,T]\times\mathbb{T}^2,~~~~
(\r,P(\r))(t,x)\in C([0,T]; W^{2,q}(\mathbb{T}^2)),\\
u\in C([0,T];H^2(\mathbb{T}^2))\cap L^2(0,T;H^3(\mathbb{T}^2)),~~
\sqrt t u\in L^\i(0,T; H^{3}(\mathbb{T}^2)),\\t u\in L^\i(0,T; W^{3,q}(\mathbb{T}^2)),~~
u_t\in L^2(0,T;H^1(\mathbb{T}^2))\\ \sqrt tu_t\in L^2(0,T; H^2(\mathbb{T}^2))\cap
L^\i(0,T;H^1(\mathbb{T}^2)),~~
t u_t\in L^\i(0,T; H^2(\mathbb{T}^2)),\\ \sqrt t\sqrt\r u_{tt}\in L^2(0,T;L^2(\mathbb{T}^2)),~~
 t\sqrt \r u_{tt}\in
L^\i(0,T;L^2(\mathbb{T}^2)),~~t\nabla u_{tt}\in
L^2(0,T;L^2(\mathbb{T}^2)).
\end{array}
\end{equation}
\end{Theorem}

\begin{Remark}
From the regularity of the solution $(\r,u)(t,x)$, it can be shown
that $(\r,u)$ is a classical solution of the system \eqref{CNS} in
$[0,T]\times\mathbb{T}^2$ (see the details in Section 5).
\end{Remark}

\begin{Remark}
If the initial data contains vacuum, then it is natural to impose the compatibility \eqref{cc} as the case of constant viscosity coefficients in \cite{CK1}.
\end{Remark}

\begin{Remark}
In Theorem \ref{theorem2}, it is not clear whether or not $u_{tt}\in
L^2(0,T;L^2(\mathbb{T}^2))$ even though one has the regularity $ t\nabla u_{tt}\in
L^2(0,T;L^2(\mathbb{T}^2))$.
\end{Remark}

\begin{Remark}
It is open to get the similar theory to the Cauchy problem or the Dirichlet problem to the 2D compressible Navier-Stokes equations \eqref{CNS}.
\end{Remark}
If the initial values are much more regular, based on Theorem
\ref{theorem2}, we can prove

\begin{Theorem}\label{theorem}
If the initial values $(\r_0,u_0)(x)$ satisfy that
\begin{equation}\label{in-d}
0\leq(\r_0(x), P(\r_0)(x))\in
H^3(\mathbb{T}^2)\times H^{3}(\mathbb{T}^2),\quad u_0(x)\in H^3(\mathbb{T}^2),\quad \int_{\mathbb{T}^2}\r_0(x) dx>0
\end{equation}
and the compatibility condition \eqref{cc}, then there exists a
unique global classical solution $(\r,u)(t,x)$ to the compressible
Navier-Stokes equations \eqref{CNS}-\eqref{initial-v} satisfying all
the properties listed in \eqref{Jan-16-1} in Theorem \ref{theorem2}
with any $2<q<\i$. Furthermore, it holds that
\begin{equation}\label{h-regu}
\begin{array}{ll}
\di u\in L^2(0,T;H^4(\mathbb{T}^2)),~~(\r,P(\r))\in C([0,T]; H^3(\mathbb{T}^2)),\\
\di \r u\in C([0,T];H^3(\mathbb{T}^2)),~~\sqrt\r \nabla^3 u\in C([0,T]; L^2(\mathbb{T}^2)).
\end{array}
\end{equation}
\end{Theorem}
\begin{Remark}
In fact, the conditions on the initial velocity $u_0$ can be weakened to $u_0\in H^2(\mathbb{T}^2)$ and $\sqrt{\r_0}\nabla^3 u_0\in L^2(\mathbb{T}^2)$ to get \eqref{h-regu}.
\end{Remark}

\begin{Remark}
In Theorem \ref{theorem}, it is not clear whether or not $u\in
C([0,T]; H^3(\mathbb{T}^2))$ even though one has $\r u\in
C([0,T]; H^3(\mathbb{T}^2))$.
\end{Remark}

\begin{Remark}
It is noted that in Theorem \ref{theorem}, the compatibility
condition \eqref{cc} is exactly same as in Theorem \ref{theorem2}.

\end{Remark}


 \vskip 2mm
\noindent\emph{Notations.} Throughout this paper, positive generic
constants are denoted by $c$ and $C$, which are independent of $\d$,
$m$ and $t\in[0,T]$, without confusion, and $C(\cdot)$ stands for
some generic constant(s) depending only on the quantity listed in
the parenthesis. For function spaces, $L^{p}(\mathbb{T}^2), 1\leq
p\leq \infty$, denote the usual Lebesgue spaces on $\mathbb{T}^2$
and $\|\cdot\|_p$ denotes its $L^p$ norm. $W^{k,p}(\mathbb{T}^2)$
denotes the $k^{th}$ order Sobolev space and
$H^{k}(\mathbb{T}^2):=W^{k,2}(\mathbb{T}^2)$.

\section{Preliminaries}

As in \cite{Kazhikhov}, we introduce the following variables. First denote the effective viscous flux by
$$
F=(2\mu+\l(\r)){\rm div}u-P(\r),
$$
and the vorticity by
$$
\o=\partial_{x_1}u_2-\partial_{x_2}u_1.
$$
Also, we define that
$$
H=\frac{1}{\r}(\mu\o_{x_1}+F_{x_2}),\qquad\qquad
L=\frac{1}{\r}(-\mu\o_{x_2}+F_{x_1}).
$$
Then the momentum equation $\eqref{CNS}_2$ can be rewritten as
\begin{equation}\label{ns1}
\left\{
\begin{array}{ll}
u_{1t}+u\cdot \nabla u_1=\frac{1}{\r}(-\mu\o_{x_2}+F_{x_1})=L,\\
u_{2t}+u\cdot \nabla u_2=\frac{1}{\r}(\mu\o_{x_1}+F_{x_2})=H.
\end{array}
\right.
\end{equation}
Then the effective viscous flux $F$ and the vorticity $\o$ solve the following system:
\begin{equation}
\left\{
\begin{array}{ll}
\o_{t}+u\cdot \nabla \o+\o{\rm div}u=H_{x_1}-L_{x_2},\\
(\f{F+P(\r)}{2\mu+\l(\r)})_{t}+u\cdot \nabla
(\f{F+P(\r)}{2\mu+\l(\r)})+(u_{1x_1})^2+2u_{1x_2}u_{2x_1}+(u_{2x_2})^2=H_{x_2}+L_{x_1}.
\end{array}
\right.
\end{equation}
Due to the continuity equation $\eqref{CNS}_1$, it holds that
\begin{equation}\label{F-omega}
\left\{
\begin{array}{ll}
\o_{t}+u\cdot \nabla \o+\o{\rm div}u=H_{x_1}-L_{x_2},\\
F_{t}+u\cdot \nabla
F-\r(2\mu+\l(\r))[F(\f{1}{2\mu+\l(\r)})^\prime+(\f{P(\r)}{2\mu+\l(\r)})^\prime]{\rm
div}u\\
\qquad+(2\mu+\l(\r))[(u_{1x_1})^2+2u_{1x_2}u_{2x_1}+(u_{2x_2})^2]=(2\mu+\l(\r))(H_{x_2}+L_{x_1}).
\end{array}
\right.
\end{equation}
Furthermore,  the system for $(H,L)$ can be derived as
\begin{equation}\label{H-L}
\left\{
\begin{array}{ll}
\r H_{t}+\r u\cdot \nabla H-\r H{\rm div}u+u_{x_2}\cdot\nabla F+\mu
u_{x_1}\cdot\nabla\o+\mu(\o{\rm
div}u)_{x_1}\\
\qquad
-\big\{\r(2\mu+\l(\r))[F(\f{1}{2\mu+\l(\r)})^\prime+(\f{P(\r)}{2\mu+\l(\r)})^\prime]{\rm
div}u\big\}_{x_2}\\
\qquad+\big\{(2\mu+\l(\r))[(u_{1x_1})^2+2u_{1x_2}u_{2x_1}+(u_{2x_2})^2]\big\}_{x_2}\\
\qquad=[(2\mu+\l(\r))(H_{x_2}+L_{x_1})]_{x_2}+\mu(H_{x_1}-L_{x_2})_{x_1},\\
\r L_{t}+\r u\cdot \nabla L-\r L{\rm div}u+u_{x_1}\cdot\nabla F-\mu
u_{x_2}\cdot\nabla\o-\mu(\o{\rm
div}u)_{x_2}\\
\qquad
-\big\{\r(2\mu+\l(\r))[F(\f{1}{2\mu+\l(\r)})^\prime+(\f{P(\r)}{2\mu+\l(\r)})^\prime]{\rm
div}u\big\}_{x_1}\\
\qquad+\big\{(2\mu+\l(\r))[(u_{1x_1})^2+2u_{1x_2}u_{2x_1}+(u_{2x_2})^2]\big\}_{x_1}\\
\qquad=[(2\mu+\l(\r))(H_{x_2}+L_{x_1})]_{x_1}-\mu(H_{x_1}-L_{x_2})_{x_2}.
\end{array}
\right.
\end{equation}
In the following, we will utilize the above systems in different
steps. Note that these systems are equivalent to each other for the
smooth solution to the original system \eqref{CNS}.

Several elementary Lemmas are needed later. The first one is the Gagliardo-Nirenberg inequality which can be found in \cite{NS}.

\begin{Lemma}\label{lemma1}
$\forall h\in W_0^{1,m}(\mathbb{T}^2)$ or $h\in
W^{1,m}(\mathbb{T}^2)$ with $\di\int_{\mathbb{T}^2} hdx=0$, it holds that
\begin{equation}
\|h\|_q\leq C\|\nabla h\|_m^\a\|h\|_r^{1-\a},
\end{equation}
where $\a=(\f1r-\f1q)(\f1r-\f1m+\f12)^{-1}$, and if $m<2,$ then $q$
is between $r$ and $\f{2m}{2-m}$, that is, $q\in[r,\f{2m}{2-m}]$ if
$r<\f{2m}{2-m}$, $q\in[\f{2m}{2-m},r]$ if $r\geq\f{2m}{2-m},$ if
$m=2,$ then $q\in[r,+\i)$, if $m>2$, then $q\in[r,+\i].$
Consequently, $\forall h\in W^{1,m}(\mathbb{T}^2)$, one has
\begin{equation}
\|h\|_q\leq C(\|h\|_1+\|\nabla h\|_m^\a\|h\|_r^{1-\a}),
\end{equation}
\end{Lemma}
where $C$ is a constant which may depend on $q$.

 The following Lemma
is the Poicare inequality.
\begin{Lemma}\label{lemma2}
$\forall h\in W_0^{1,m}(\mathbb{T}^2)$ or $h\in
W^{1,m}(\mathbb{T}^2)$ with $\di\int_{\mathbb{T}^2} hdx=0$,  if
$1\leq m<2,$ then
\begin{equation}
\|h\|_{\f{2m}{2-m}}\leq C(2-m)^{-\f12}\|\nabla h\|_m,
\end{equation}
where the positive constant $C$ is independent of $m.$
\end{Lemma}
The following Lemma follows from Lemma \ref{lemma2}, of which proof
can be found in \cite{Kazhikhov}.
\begin{Lemma}\label{lemma3}
$\forall h\in W^{1,\f{2m}{m+\eta}}(\mathbb{T}^2)$ with $m\geq2$ and
$0<\eta\leq1$, we have
\begin{equation}
\|h\|_{2m}\leq C(\|h\|_1+m^{\f12}\| h\|_{2(1-\v)}^s\|\nabla
h\|_{\f{2m}{m+\eta}}^{1-s}),
\end{equation}
where $\v\in[0,\f12], s=\f{(1-\v)(1-\eta)}{m-\eta(1-\v)}$  and the
positive constant $C$ is independent of $m.$
\end{Lemma}

\section{Approximate solutions}
\setcounter{equation}{0}

In this section, we construct a sequence of approximate solutions by making use of the theory
of Vaigant-Kazhikhov \cite {Kazhikhov} and derive some uniform a-priori estimates which are
necessary to prove Theorem \ref{theorem2}. To this end, we need a careful approximation of the
initial data.

\underline{Step 1. Approximation of initial data:}
To apply the theory of
Vaigant-Kazhikhov \cite{Kazhikhov}, we approximate
of the initial data in \eqref{in-d} as follows. First. the initial density and pressure can be approximated as
\begin{equation}\label{app-in-d}
\r_0^\d=\r_0+\d,  \qquad P_0^\d=P(\r_0)+\d,
\end{equation}
for any small positive constant $\d>0$.
To approximate the initial velocity, we define $u_0^\d$ to be the unique solution to the following elliptic problem
\begin{equation}\label{new-cc}
  \mathcal{L}_{\r_0^\d}u_0^\d=\nabla P_0^\d+\sqrt{\r_0}g
\end{equation}
with the periodic boundary conditions on $\mathbb{T}^2$ and $\di \int_{\mathbb{T}^2} u_0^\d dx=\int_{\mathbb{T}^2} u_0dx:=\bar u_0.$ It should be noted that $u_0^\d$ is uniquely determined due to the compatibility condition \eqref{cc}.

It follows from \eqref{new-cc} that
\begin{equation}\label{new-cc1}
  \mathcal{L}_{\r_0} u_0^\d=-\nabla\big[(\l(\r_0^\d)-\l(\r_0)){\rm div} u_0^\d\big]+\nabla P_0^\d+\sqrt{\r_0} g.
\end{equation}
By the elliptic regularity, it holds that
\begin{equation}\label{ud1}
\begin{array}{ll}
\di \|u_0^\d-\bar u_0\|_{H^2(\mathbb{T}^2)}\\
\di\leq C\Big[\|\l(\r_0^\d)-\l(\r_0)\|_\i\|\nabla({\rm div} u_0^\d)\|_2+\|\nabla(\l(\r_0^\d)-\l(\r_0))\|_\i\|{\rm div} u_0^\d\|_2+\|\nabla P_0^\d\|_2+\|\sqrt{\r_0} g\|_2\Big]\\
\di\leq C\Big[\d\|u_0^\d\|_{H^2(\mathbb{T}^2)}+\|P_0\|_{H^1(\mathbb{T}^2)}+\|\sqrt{\r_0}\|_{L^\i(\mathbb{T}^2)} \|g\|_2\Big]\\
\leq C\Big[\d\| u_0^\d\|_{H^2(\mathbb{T}^2)}+1\Big].
\end{array}
\end{equation}
where the generic positive constant $C$ is independent of $\d>0.$

Therefore, if $\d\ll1$, then \eqref{ud1} yields that
\begin{equation}\label{h2}
\|u_0^\d\|_{H^2(\mathbb{T}^2)}\di\leq C
\end{equation}
where the positive constant $C$ is independent of $0<\d\ll1.$

Due to the compatibility condition \eqref{cc} and \eqref{new-cc}, it holds that
\begin{equation}\label{h5}
\begin{array}{ll}
\di \mathcal{L}_{\r_0}(u_0^\d-u_0)
&\di=-\nabla\big[(\l(\r_0^\d)-\l(\r_0)){\rm div} u_0^\d\big]:=\T^\d.
\end{array}
\end{equation}
Therefore, by the elliptic regularity, \eqref{app-in-d} and \eqref{h2}, one can get that
\begin{equation}\label{h6}
\begin{array}{ll}
\|u_0^\d-u_0\|_{H^2(\mathbb{T}^2)}\leq C\|\T^\d\|_2\\
\qquad\di\leq\di C\Big[\|\l(\r_0^\d)-\l(\r_0)\|_{L^\i(\mathbb{T}^2)}\| \nabla^2u_0^\d\|_{2}+\|\nabla(\l(\r_0^\d)-\l(\r_0))\|_{L^\i(\mathbb{T}^2)}\| {\rm div} u_0^\d\|_{2}\Big]\\
\qquad\di\leq C\Big[\|\l(\r_0^\d)-\l(\r_0)\|_{L^\i(\mathbb{T}^2)}+\|\nabla(\l(\r_0^\d)-\l(\r_0))\|_{L^\i(\mathbb{T}^2)}\Big]\\ \qquad\di\leq C\d~\rightarrow 0, \qquad {\rm as}~~\d\rightarrow 0.
\end{array}
\end{equation}
For the initial data $(\r_0^\d, P_0^\d, u_0^\d)$ constructed above
for each fixed $\d>0$, it is proved in \cite{Kazhikhov} that the
compressible Navier-Stokes equations \eqref{CNS} with $\b>3$ has a
unique global strong solution $(\r^\d,u^\d)$ such that
$c_{\d}\leq\r^\d\leq C_\d$ for some positive constants $c_\d, C_\d$
depending on $\d$. In the following, we will derive the uniform
bound to $(\r^\d,u^\d)$ with respect to $\d$ and then pass the limit
$\d\rightarrow0$ to get the classical solution which may contain
vacuum states in an open set of $\mathbb{T}^2$. It should be noted
that in comparison with estimates presented in \cite{Kazhikhov}, we
will obtain uniform estimates with respective to the lower bound of
the density such that vacuum is permitted in these estimates. To
this end, the compatibility condition \eqref{cc} will be crucial.

For simplicity of notations, we will omit the superscript $^\d$ of
$(\r^\d,u^\d)$ in the following in the case of no confusions.

\underline{Step 2. Elementary energy estimates:}
\begin{Lemma}\label{lemma-ee}
There exists a positive constant $C$ depending on $(\r_0,u_0)$, such that
\begin{equation}
\sup_{t\in[0,T]}\big(\|\sqrt\r u\|^2_2+\|\r\|^\g_\g\big)+\int_0^T\big(\|\nabla u\|_2^2+\|\o\|_2^2+\|(2\mu+\l(\r))^\f12{\rm div} u\|_2^2\big)dt\leq C.
\end{equation}
\end{Lemma}
{\bf Proof:} Multiplying the equation $\eqref{ns1}_i$ by $\r u_i, (i=1,2)$, summing the resulting equations and then integrating over $\mathbb{T}^2$ and using the continuity equation $\eqref{CNS}_1$, it holds that
$$
\f{d}{dt}\int \r|u|^2dx+ \int(\mu\o^2+(2\mu+\l(\r))({\rm div} u)^2)dx+\int u\cdot\nabla Pdx=0.
$$
Multiplying the continuity equation $\eqref{CNS}_1$ by $\f{1}{\g-1}\r^{\g-1}$ and then integrating over $\mathbb{T}^2$ yields that
$$
\f{d}{dt}\int\f{\r^\g}{\g-1}dx+\int P{\rm div}udx=0.
$$
Therefore, combining the above two estimates and then integrating over $[0,t]$ with respect to $t$, we obtain
\begin{equation}\label{Eee}
\begin{array}{ll}
  \di\int(\f12\r |u|^2+\f{1}{\g-1}\r^\g)dx+\int_0^T\int
(\mu\o^2+(2\mu+\l(\r))({\rm div} u)^2)dxdt\\
\di=\int(\f12\r_0^\d |u_0^\d|^2+\f{1}{\g-1}(\r_0^\d)^\g)dx\\
\di \leq C\Big[\|\r_0^\d\|_{H^3(\mathbb{T}^2)}\|u_0^\d\|^2_{H^2(\mathbb{T}^2)}+\|\r_0^\d\|_{H^3(\mathbb{T}^2)}^\g\Big]\leq C.
\end{array}
\end{equation}
Denote
\begin{equation}\label{Jan-18-1}
\phi(t)=\int (\mu\o^2+(2\mu+\l(\r))({\rm div} u)^2)dx,\qquad
t\in[0,T].
\end{equation}
Then
$$
\|\nabla u\|_{2}^2(t)\leq C\big[\|\o\|_2^2(t)+\|{\rm div} u\|_2^2(t)\big] \leq C\phi(t)\in L^1(0,T).
$$
Thus the proof of Lemma \ref{lemma-ee} is completed. $\hfill\Box$

\underline{Step 3. Density estimates:} Applying the operator $div$
to the momentum equation $\eqref{CNS}_2$, we have
\begin{equation}\label{vis-f}
[{\rm div}(\r u)]_t+{\rm div}[{\rm div}(\r u\otimes u)]=\Delta F.
\end{equation}
Consider the following two elliptic problems:
\begin{equation}\label{xi}
\Delta \x={\rm div}(\r u),\qquad \int\x dx=0,
\end{equation}
\begin{equation}\label{eta}
\Delta \eta={\rm div}[{\rm div}(\r u\otimes u)],\qquad \int\eta
dx=0,
\end{equation}
both with the periodic boundary condition on the torus
$\mathbb{T}^2.$

By the elliptic estimates and H${\rm\ddot{o}}$lder inequality, it
holds that
\begin{Lemma}\label{lemma4}
\begin{itemize}
\item[(1)] $\|\nabla\x\|_{2m}\leq Cm\|\r\|_{\f{2mk}{k-1}}\|u\|_{2mk},$
for any $k>1,m\geq1;$
\item[(2)] $\|\nabla\x\|_{2-r}\leq C\|\sqrt\r u\|_{2}\|\r\|^{\f12}_{\f{2-r}{r}},$
for any $0< r<1;$
\item[(3)] $\|\eta\|_{2m}\leq Cm\|\r\|_{\f{2mk}{k-1}}\|u\|^2_{4mk},$
for any $k>1,m\geq1;$
\end{itemize} where $C$ are positive constants independent of $m,k$
and $r$.
\end{Lemma}
{\bf Proof:} (1) By the elliptic estimates to the equation
\eqref{xi} and then using the H${\rm\ddot{o}}$lder inequality, we
have for any $k>1,m\geq1$,
$$
\|\nabla\x\|_{2m}\leq C m\|\r u\|_{2m}\leq
Cm\|\r\|_{\f{2mk}{k-1}}\|u\|_{2mk}.
$$
Similarly, the statements (2) and (3) can be proved. $\hfill\Box$

Based on Lemmas \ref{lemma1}-\ref{lemma3} and Lemma \ref{lemma4}, it holds that

\begin{Lemma}\label{lemma5}
\begin{itemize}
\item[(1)] $\|\x\|_{2m}\leq Cm^{\f12}\|\nabla\x\|_{\f{2m}{m+1}}\leq Cm^{\f12}\|\r\|_m^{\f12},$
for any $m\geq2;$
\item[(2)] $\|u\|_{2m}\leq C\left[m^{\f12}\|\nabla u\|_{2}+1\right],$
for any $m\geq2;$
\item[(3)] $\|\nabla\x\|_{2m}\leq C\left[m^{\f32}k^{\f12}\|\r\|_{\f{2mk}{k-1}}\phi(t)^{\f12}+m\|\r\|_{\f{2mk}{k-1}}\right],$
for any $k>1,m\geq1;$
\item[(4)] $\|\eta\|_{2m}\leq C\left[m^2k\|\r\|_{\f{2mk}{k-1}}\phi(t)+m\|\r\|_{\f{2mk}{k-1}}\right],$
for any $k>1,m\geq1;$
\end{itemize} where $C$ are positive constants independent of $m,k$.
\end{Lemma}
{\bf Proof:} (1) By Lemma \ref{lemma2} and Lemma \ref{lemma4} (2),
it holds that
$$
\begin{array}{ll}
\di \|\x\|_{2m}\leq Cm^{\f12}\|\nabla\x\|_{\f{2m}{m+1}}&\di\leq
Cm^{\f12}\|\sqrt\r u\|_2\|\r\|_m^{\f12}\\
&\di\leq Cm^{\f12}\|\r\|_m^{\f12},
\end{array}
$$
where in the last inequality one has used the elementary energy
estimates \eqref{Eee}.

(2). From the conservative form of the compressible Navier-Stokes
equations \eqref{CNS} and the periodic boundary conditions, we have
$$
\f{d}{dt}\int\r(t,x) dx=\f{d}{dt}\int\r u(t,x) dx=0,
$$
that is,
$$
\int\r(t,x)dx=\int\r_0(x) dx, \qquad \int\r u(t,x)dx=\int\r_0u_0(x)
dx, \qquad \forall t\in[0,T].
$$
By Lemma \ref{lemma2}, it follows that
\begin{equation}\label{u1}
\|u\|_{2m}\leq \|u-\bar u\|_{2m}+\|\bar u\|_{2m}\leq
Cm^{\f12}\|\nabla u\|_{\f{2m}{m+1}}+|\bar u|,
\end{equation}
where $m>2$ and $\bar u=\bar u(t)=\di\int u(t,x)dx.$

On the other hand, we have
\begin{equation}\label{bar-u1}
|\int\r(u-\bar u) dx|\leq \|\r\|_{\g}\|u-\bar u\|_{\f{\g}{\g-1}}\leq
C\|\nabla u\|_2,
\end{equation}
where in the last inequality we have used the elementary energy
estimates \eqref{Eee} and the Poincare inequality.

Note that
\begin{equation}\label{bar-u2}
|\int\r(u-\bar u) dx|=|\int\r_0 u_0dx-\bar u\int\r_0(x)dx|\geq |\bar
u|\int\r_0dx-|\int\r_0 u_0dx|.
\end{equation}
Combining \eqref{bar-u1} with \eqref{bar-u2} implies that
\begin{equation}\label{bar-u3}
|\bar u|\leq \f{|\int\r_0 u_0dx|}{\int\r_0dx}+
\f{C\|\nabla u\|_2}{\int\r_0dx}.
\end{equation}
Substituting \eqref{bar-u3} into \eqref{u1} completes the proof of
Lemma \ref{lemma5} (2).

The assertions (3) and (4) in Lemma \ref{lemma5} are
direct consequences of Lemma \ref{lemma5} (2) and Lemma
\ref{lemma4} (1), (3), respectively. Thus the proof of Lemma \ref{lemma5} is completed. $\hfill\Box$

Substituting \eqref{xi} and \eqref{eta} into \eqref{vis-f} yields
that
\begin{equation}
\Delta\Big(\x_t+\eta-F+\int F(t,x)dx\Big)=0.
\end{equation}
Thus, it holds that
\begin{equation}
\x_t+\eta-F+\int F(t,x)dx=0.
\end{equation}
It follows from the definition of the effective viscous flux $F$ that
\begin{equation}
\x_t-(2\mu+\l(\r)){\rm div}u+P(\r)+\eta+\int F(t,x)dx=0.
\end{equation}
Then the continuity equation $\eqref{CNS}_1$ yields that
\begin{equation}
\x_t+\f{2\mu+\l(\r)}{\r}(\r_t+u\cdot\nabla \r)+P(\r)+\eta+\int
F(t,x)dx=0.
\end{equation}
Define
\begin{equation}\label{theta}
\theta(\r)=\int_1^\r\f{2\mu+\l(s)}{s}ds=2\mu\ln\r+\f{1}{\b}(\r^\b-1).
\end{equation}
Then we obtain the following transport equation
\begin{equation}\label{transport-e}
(\x+\t(\r))_t+u\cdot\nabla(\x+\t(\r))+P(\r)+\eta-u\cdot\nabla\x+\int
F(t,x)dx=0.
\end{equation}
\begin{Lemma}\label{lemma-rho}
  For any $k\geq1,$ it holds that
  \begin{equation}\label{density-e}
    \sup_{t\in[0,T]}\|\r(t,\cdot)\|_k\leq C k^{\f{2}{\b-1}}.
  \end{equation}
\end{Lemma}
{\bf Proof:} Multiplying the equation \eqref{transport-e} by
$\r[(\x+\t(\r))_+]^{2m-1}$ with $m\geq4$ being integer, here and in what follows, the
notation $(\cdots)_+$ denotes the positive part of $(\cdots)$, one can get that
\begin{equation}\label{trans-e1}
\begin{array}{ll}
\di \f{1}{2m}\f{d}{dt}\int\r[(\x+\t(\r))_+]^{2m}dx+\int\r
P(\r)[(\x+\t(\r))_+]^{2m-1}dx=-\int\r \eta[(\x+\t(\r))_+]^{2m-1}dx\\
\di  +\int\r u\cdot\nabla \x[(\x+\t(\r))_+]^{2m-1}dx-\int F(t,x)dx
\int \r [(\x+\t(\r))_+]^{2m-1}dx.
\end{array}
\end{equation}
Denote
\begin{equation}
f(t)=\big\{\int\r[(\x+\t(\r))_+]^{2m} dx\big\}^{\f{1}{2m}},\qquad
t\in[0,T].
\end{equation}
Now we estimate the terms on the right hand side of
\eqref{trans-e1}. First,
\begin{equation}\label{rho-1}
\begin{array}{ll}
\di |-\int\r \eta[(\x+\t(\r))_+]^{2m-1}dx|\leq
\int\r^{\f{1}{2m}}|\eta|\big[\r(\x+\t(\r))^{2m}_+\big]^{\f{2m-1}{2m}}dx\\
\qquad\di\leq
\|\r\|_{2m\b+1}^{\f{1}{2m}}\|\eta\|_{2m+\f{1}{\b}}\|\r(\x+\t(\r))^{2m}_+\|_1^{\f{2m-1}{2m}}\\
\qquad\di\leq
C\|\r\|_{2m\b+1}^{\f{1}{2m}}\Big[(m+\f{1}{2\b})^2k\|\r\|_{\f{2(m+\f{1}{2\b})k}{k-1}}\phi(t)+(m+\f{1}{2\b})\|\r\|_{\f{2(m+\f{1}{2\b})k}{k-1}}\Big]f(t)^{2m-1}\\
\qquad\di\leq
C\|\r\|_{2m\b+1}^{1+\f{1}{2m}}f(t)^{2m-1}\big[m^2\phi(t)+m\big],
\end{array}
\end{equation}
where $\phi(t)$ is defined as in \eqref{Jan-18-1} and in the last
inequality we have taken $k=\f{\b}{\b-1}.$

Next, for $\f{1}{2m\b+1}+\f1p+\f1q=1$ with $p,q\geq 1$, one has
\begin{equation}\label{rho-2}
\begin{array}{ll}
\di |\int\r u\cdot\nabla\x[(\x+\t(\r))_+]^{2m-1}dx|\leq
\int\r^{\f{1}{2m}}|u||\nabla\x|\big[\r(\x+\t(\r))^{2m}_+\big]^{\f{2m-1}{2m}}dx\\
\qquad\di\leq
\|\r\|_{2m\b+1}^{\f{1}{2m}}\|u\|_{2mp}\|\nabla\x\|_{2mq}\|\r(\x+\t(\r))^{2m}_+\|_1^{\f{2m-1}{2m}}\\
\qquad\di\leq
C\|\r\|_{2m\b+1}^{\f{1}{2m}}\Big[(mp)^{\f12}\|\nabla u\|_2+1\Big]\Big[(mq)^{\f32}k^{\f12}\|\r\|_{\f{2mqk}{k-1}}\phi(t)^{\f12}+m\|\r\|_{\f{2mqk}{k-1}}\Big]f(t)^{2m-1}\\
\qquad\di\leq
C\|\r\|_{2m\b+1}^{1+\f{1}{2m}}f(t)^{2m-1}\big[m^{\f12}\phi(t)^{\f12}+1\big]\big[m^{\f32}\phi(t)^{\f12}+m\big]\\
\qquad\di\leq
C\|\r\|_{2m\b+1}^{1+\f{1}{2m}}f(t)^{2m-1}\big[m^2\phi(t)+m\big],
\end{array}
\end{equation}
where in the third inequality one has chosen $p=q=\f{2m\b+1}{m\b}$
and $k=\f{\b}{\b-1}.$

Then it follows that
\begin{equation}\label{rho-3}
\begin{array}{ll}
\di |-\int F(t,x)dx \int \r [(\x+\t(\r))_+]^{2m-1}dx|\\
\di\leq\int|(2\mu+\l(\r)){\rm div} u-P(\r)|dx\int
\r^{\f{1}{2m}}\big[\r(\x+\t(\r))^{2m}_+\big]^{\f{2m-1}{2m}}dx\\
\di\leq\Big[(\int(2\mu+\l(\r))({\rm
div}u)^2dx)^{\f12}(\int(2\mu+\l(\r))dx)^{\f12}+\int
P(\r)dx\Big]\|\r\|_1^{\f{1}{2m}}\|\r(\x+\t(\r))^{2m}_+\|_1^{\f{2m-1}{2m}}\\
\di\leq C\Big[\phi(t)^{\f12}+\phi(t)^{\f12}(\int\r^\b
dx)^{\f12}+1\Big]f(t)^{2m-1}\\
\di\leq
C\Big[\phi(t)^{\f12}+\phi(t)^{\f12}\|\r\|_{2m\b+1}^{\f\b2}+1\Big]f(t)^{2m-1}.
\end{array}
\end{equation}
Substituting \eqref{rho-1}, \eqref{rho-2} and \eqref{rho-3} into
\eqref{trans-e1} yields that
\begin{equation}
\begin{array}{ll}
\di\f{1}{2m}\f{d}{dt}(f^{2m}(t))+\int\r P(\r)[(\x+\t(\r))_+]^{2m-1}dx\\
\di\leq
C\|\r\|_{2m\b+1}^{1+\f{1}{2m}}f(t)^{2m-1}\big[m^2\phi(t)+m\big]+C\Big[\phi(t)^{\f12}+\phi(t)^{\f12}\|\r\|_{2m\b+1}^{\f\b2}+1\Big]f(t)^{2m-1}.
\end{array}
\end{equation}
Then it holds that
\begin{equation}
\f{d}{dt}f(t)\leq
C\Big[1+\phi(t)^{\f12}+\phi(t)^{\f12}\|\r\|_{2m\b+1}^{\f\b2}+\big(m^2\phi(t)+m\big)\|\r\|_{2m\b+1}^{1+\f{1}{2m}}\Big].
\end{equation}
Integrating the above inequality over $[0,t]$ gives that
\begin{equation}\label{f}
f(t)\leq
f(0)+C\Big[1+\int_0^t\phi(\tau)^{\f12}\|\r\|_{2m\b+1}^{\f\b2}(\tau)d\tau+\int_0^t\big(m^2\phi(\tau)+m\big)\|\r\|_{2m\b+1}^{1+\f{1}{2m}}(\tau)d\tau\Big].
\end{equation}
Now we calculate the quantity
$$f(0)=\Big(\int\r^\d_0[(\x^\d_0+\t(\r^\d_0))_+]^{2m} dx\Big)^{\f{1}{2m}}.$$
By Lemma \ref{lemma4} (1) with $t=0$, we can easily get
$$
\|\x^\d_0\|_{L^\i}\leq C.
$$
Furthermore, by the definition of
$\t(\r^\d_0)=2\mu\ln\r^\d_0+\f1\b((\r^\d_0)^\b-1)$, we have
$$
\x^\d_0+\t(\r^\d_0)\rightarrow-\i, \quad{\rm as}\quad \r^\d_0\rightarrow0+.
$$
Thus there exists a positive constant $\s$, such that if
$0\leq\r^\d_0\leq\s$, then
$$
(\x^\d_0+\t(\r^\d_0))_+\equiv0.
$$
Now one has
\begin{equation}\label{f0}
\begin{array}{ll}
f(0)&\di =\Big[\Big(\int_{[0\leq\r_0\leq
\s]}+\int_{[\s\leq\r^\d_0\leq M]}\Big)\r^\d_0(\x^\d_0+\t(\r^\d_0))_+^{2m}
dx\Big]^{\f{1}{2m}}\\
&\di=\Big[\int_{[\s\leq\r^\d_0\leq M]}\r^\d_0(\x^\d_0+\t(\r^\d_0))_+^{2m}
dx\Big]^{\f{1}{2m}}\leq C(\s,M),
\end{array}
\end{equation}
where the positive constant $C(\s,M)$ is independent of $\d$ and $m$.

It follows from \eqref{f} and \eqref{f0} that
\begin{equation}\label{f1}
f(t)\leq
C\Big[1+\int_0^t\phi(\tau)^{\f12}\|\r\|_{2m\b+1}^{\f\b2}(\tau)d\tau+\int_0^t\big(m^2\phi(\tau)+m\big)\|\r\|_{2m\b+1}^{1+\f{1}{2m}}(\tau)d\tau\Big].
\end{equation}
Set $\Omega_1(t)=\{x\in\mathbb{T}^2|\r(t,x)>2\}$ and
$\Omega_2(t)=\{x\in\Omega_1(t)|\xi(t,x)+\t(\r)(t,x)>0\}$. Then one
has
\begin{equation}
\begin{array}{ll}
&\di\|\r\|_{2m\b+1}^\b(t)=\Big(\int\r^{2m\b+1}dx\Big)^{\f{\b}{2m\b+1}}=\Big(\int_{\Omega_1(t)}\r^{2m\b+1}dx+\int_{\mathbb{T}^2\setminus\Omega_1(t)}\r^{2m\b+1}dx\Big)^{\f{\b}{2m\b+1}}\\
&\di~~\leq\Big(\int_{\Omega_1(t)}\r^{2m\b+1}dx\Big)^{\f{\b}{2m\b+1}}+C\leq
C\Big(\int_{\Omega_1(t)}\r|\t(\r)|^{2m}dx\Big)^{\f{\b}{2m\b+1}}+C\\
&\di~~=
C\Big(\int_{\Omega_2(t)}\r|\t(\r)+\x-\x|^{2m}dx+\int_{\Omega_1(t)\setminus\Omega_2(t)}\r|\t(\r)|^{2m}dx\Big)^{\f{\b}{2m\b+1}}+C\\
&\di~~\leq
C\Big(\int_{\Omega_2(t)}\r(\t(\r)+\x)^{2m}dx+\int_{\Omega_2(t)}\r|\x|^{2m}dx+\int_{\Omega_1(t)\setminus\Omega_2(t)}\r|\x|^{2m}dx\Big)^{\f{\b}{2m\b+1}}+C\\
&\di~~\leq
C\Big(f(t)^{2m}+\int_{\mathbb{T}^2}\r|\x|^{2m}dx\Big)^{\f{\b}{2m\b+1}}+C\leq
C\Big[f(t)+\Big(\int_{\mathbb{T}^2}\r|\x|^{2m}dx\Big)^{\f{\b}{2m\b+1}}+1\Big].
\end{array}
\end{equation}
Note that
\begin{equation}
\begin{array}{ll}
\di\Big(\int_{\mathbb{T}^2}\r|\x|^{2m}dx\Big)^{\f{\b}{2m\b+1}}&\di
\leq
\|\r\|_{2m\b+1}^{\f{\b}{2m\b+1}}\||\x|^{2m}\|^{\f{\b}{2m\b+1}}_{\f{2m\b+1}{2m\b}}=\|\r\|_{2m\b+1}^{\f{\b}{2m\b+1}}\|\x\|^{\f{2m\b}{2m\b+1}}_{2m+\f{1}{\b}}\\
&\di \leq\|\r\|_{2m\b+1}^{\f{\b}{2m\b+1}}\Big[C(m+\f{1}{2\b})^{\f12}\|\r\|_{m+\f{1}{2\b}}^{\f12}\Big]^{\f{2m\b}{2m\b+1}}
\leq Cm^{\f12}\|\r\|_{2m\b+1}^{\f{\b(m+1)}{2m\b+1}},
\end{array}
\end{equation}
Then one can get
\begin{equation}
\begin{array}{ll}
\|\r\|_{2m\b+1}^\b(t)&\di\leq
C\Big[1+f(t)+m^{\f12}\|\r\|_{2m\b+1}^{\f{\b(m+1)}{2m\b+1}}(t)\Big]\\
&\di \leq
\f12\|\r\|_{2m\b+1}^\b(t)+C\Big(1+f(t)+m^{\f{m\b+\f12}{m(2\b-1)}}\Big).
\end{array}
\end{equation}
Thus it holds that
\begin{equation}
\begin{array}{ll}
\|\r\|_{2m\b+1}^\b(t)&\di\leq C\Big[f(t)+m^{\f{\b}{2\b-1}}\Big]\\
&\di\leq
C\Big[m^{\f{\b}{2\b-1}}+\int_0^t\phi(\tau)^{\f12}\|\r\|_{2m\b+1}^{\f\b2}(\tau)d\tau+\int_0^t\big(m^2\phi(\tau)+m\big)\|\r\|_{2m\b+1}^{1+\f{1}{2m}}(\tau)d\tau\Big]\\
&\di\leq
C\Big[m^{\f{\b}{2\b-1}}+\int_0^t\|\r\|_{2m\b+1}^{\b}(\tau)d\tau+\int_0^t\big(m^2\phi(\tau)+m\big)\|\r\|_{2m\b+1}^{1+\f{1}{2m}}(\tau)d\tau\Big].
\end{array}
\end{equation}
Applying Gronwall's inequality yields that
\begin{equation}
\|\r\|_{2m\b+1}^\b(t)\leq
C\Big[m^{\f{\b}{2\b-1}}+\int_0^t\big(m^2\phi(\tau)+m\big)\|\r\|_{2m\b+1}^{1+\f{1}{2m}}(\tau)d\tau\Big].
\end{equation}
Denote
$$
y(t)=m^{-\f{2}{\b-1}}\|\r\|_{2m\b+1}(t).
$$
Then it holds that
$$
\begin{array}{ll}
\di y^\b(t)&\di \leq
C\Big[m^{\f{\b(1-3\b)}{(2\b-1)(\b-1)}}+m^{\f{1}{m(\b-1)}}\int_0^t\phi(\tau)y(\tau)^{1+\f{1}{2m}}d\tau+m^{\f{1}{m(\b-1)}-1}\int_0^ty(\tau)^{1+\f{1}{2m}}d\tau\Big]\\
&\di \leq C\Big[1+\int_0^t(\phi(\tau)+1)y^\b(\tau)d\tau\Big].
\end{array}
$$
So applying the Gronwall's inequality again yields that
$$
y(t)\leq C,\quad \forall t\in[0,T],
$$
that is,
$$
\|\r\|_{2m\b+1}(t)\leq Cm^{\f{2}{\b-1}},\quad \forall t\in[0,T].
$$
Equivalently,  \eqref{density-e} holds. Thus Lemma \ref{lemma-rho} is proved. $\hfill\Box$

\underline{Step 4: First-order derivative estimates of the
velocity.}

\begin{Lemma}\label{lemma-u-der}
  There exists a positive constant $C$, such that
  \begin{equation}
    \sup_{t\in[0,T]}\int(\mu\omega^2+\f{F^2}{2\mu+\l(\r)}) dx+\int_0^T\int\r(H^2+L^2)
     dxdt\leq C.
  \end{equation}
\end{Lemma}
{\bf Proof:} Multiplying the equation $\eqref{F-omega}_1$ by $\mu\omega$, the
equation $\eqref{F-omega}_2$ by $\f{F}{2\mu+\l(\r)}$, respectively,
and then summing the resulted equations together, one has
\begin{equation}\label{ue1}
\begin{array}{ll}
\di\f12\f{d}{dt}\int(\mu\omega^2+\f{F^2}{2\mu+\l(\r)}) dx+\f{\mu}{2}\int\omega^2{\rm div }u dx-\f12\int\r F^2(\f{1}{2\mu+\l(\r)})^\prime{\rm div }u dx \\
\di \quad -\f12\int F^2\f{{\rm div}u}{2\mu+\l(\r)} dx-\int\r F({\rm
div} u)(\f{P(\r)}{2\mu+\l(\r)})^\prime dx+\int
F[(u_{1x_1})^2+2u_{1x_2}u_{2x_1}+(u_{2x_2})^2]dx\\
\di\quad=-\int\r(H^2+L^2) dx.
\end{array}
\end{equation}
Notice that
$$
\begin{array}{ll}
\di (u_{1x_1})^2+2u_{1x_2}u_{2x_1}+(u_{2x_2})^2&\di =(u_{1x_1}+u_{2x_2})^2+2(u_{1x_2}u_{2x_1}-u_{1x_1}u_{2x_2})\\
&\di=({\rm div} u)^2+2(u_{1x_2}u_{2x_1}-u_{1x_1}u_{2x_2})\\
&\di=({\rm div}
u)\left(\f{F}{2\mu+\l(\r)}+\f{P(\r)}{2\mu+\l(\r)}\right)+2(u_{1x_2}u_{2x_1}-u_{1x_1}u_{2x_2}),
\end{array}
$$
then one has
\begin{equation}\label{ue1}
\begin{array}{ll}
\di\f12\f{d}{dt}\int(\mu\omega^2+\f{F^2}{2\mu+\l(\r)}) dx+\int\r(H^2+L^2) dx=-\f{\mu}{2}\int\omega^2{\rm div }u dx\\
\di \quad +\f12\int F^2({\rm
div}u)\Big[\r(\f{1}{2\mu+\l(\r)})^\prime-\f{1}{2\mu+\l(\r)}\Big]
dx+\int F({\rm div}
u)\Big[\r(\f{P(\r)}{2\mu+\l(\r)})^\prime-\f{P(\r)}{2\mu+\l(\r)}\Big]
dx\\
\di \quad -\int 2F(u_{1x_2}u_{2x_1}-u_{1x_1}u_{2x_2})dx.
\end{array}
\end{equation}
Set
$$
Z^2(t)=\int(\mu\omega^2+\f{F^2}{2\mu+\l(\r)}) dx,
$$
and
$$
\varphi^2(t)=\int\r(H^2+L^2)
dx=\int\frac{1}{\r}\big[(\mu\o_{x_1}+F_{x_2})^2+(-\mu\o_{x_2}+F_{x_1})^2\big]
dx.
$$
Then it follows  that for $0<r\leq \f12,$
\begin{equation}\label{fact1}
\di \|\nabla(F,\omega)\|_{2(1-r)} \di \leq
C\varphi(t)\|\r\|_{\f{1-r}{r}}^{\f12}\leq C\varphi(t)
(\f{1-r}{r})^{\f{1}{\b-1}}\leq C\varphi(t)r^{\f{1}{1-\b}},
\end{equation}
and
\begin{equation}\label{fact2}
\begin{array}{ll}
\di \|\nabla u\|_2+\|\omega\|_2+\|{\rm div}
u\|_2+\|(2\mu+\l(\r))^{\f12}{\rm div} u\|_2\\
\di \leq C\Big[Z(t)+\Big(\int\f{P(\r)^2}{2\mu+\l(\r)}
dx\Big)^{\f12}\Big]\leq C(Z(t)+1).
\end{array}
\end{equation}

Now we estimate the four terms on the right hand side of
\eqref{ue1}. First, by the interpolation inequality and Lemma
\ref{lemma2}, \eqref{fact1} and \eqref{fact2}, for $0<\v\leq\f14$,
it holds that
\begin{equation}\label{ue2}
\begin{array}{ll}
\di |-\f{\mu}{2}\int\omega^2{\rm div }u dx|\leq C\|{\rm div}
u\|_2\|\o\|_4^2 \leq
C(Z(t)+1)\|\o\|_2^{\f{1-3\v}{1-2\v}}\|\nabla\o\|_{2(1-\v)}^{\f{1-\v}{1-2\v}}\\
\di \qquad\qquad\leq
C(Z(t)+1)Z(t)^{\f{1-3\v}{1-2\v}}\varphi(t)^{\f{1-\v}{1-2\v}}\v^{\f{1-\v}{(1-\b)(1-2\v)}}\\
\di\qquad\qquad \leq \a \varphi^2(t)+C_\a
Z(t)^2(Z(t)+1)^{\f{2(1-2\v)}{1-3\v}}\v^{\f{2}{1-\b}\f{1-\v}{1-3\v}}\\
\di\qquad\qquad \leq \a \varphi^2(t)+C_\a
(Z(t)^2+1)^{2+\f{\v}{1-3\v}}\v^{\f{2}{1-\b}\f{1-\v}{1-3\v}},
\end{array}
\end{equation}
where and in the sequel $\a>0$ is a small positive constant to be determined and
$C_\a$ is a positive constant depending on $\a$.

Next, one has
\begin{equation}\label{ue3}
\begin{array}{ll}
\di |\f12\int F^2{\rm
div}u\Big[\r(\f{1}{2\mu+\l(\r)})^\prime-\f{1}{2\mu+\l(\r)}\Big]
dx|\\
\di =|\f12\int
F^2\left(\f{F}{2\mu+\l(\r)}+\f{P(\r)}{2\mu+\l(\r)}\right)\f{2\mu+\l(\r)+\r\l^\prime(\r)}{(2\mu+\l(\r))^2}
dx|\\
\di \leq C\int|F|^2
\left(\f{|F|}{2\mu+\l(\r)}+\f{P(\r)}{2\mu+\l(\r)}\right)dx \leq C\left(1+\int \f{|F|^3}{2\mu+\l(\r)}dx\right),
\end{array}
\end{equation}
and
\begin{equation}\label{ue4}
\begin{array}{ll}
\di |\f12\int F{\rm
div}u\Big[\r(\f{P(\r)}{2\mu+\l(\r)})^\prime-\f{P(\r)}{2\mu+\l(\r)}\Big]
dx|\\
\di =|\f12\int
F\left(\f{F}{2\mu+\l(\r)}+\f{P(\r)}{2\mu+\l(\r)}\right)\f{P(\r)(2\mu+\l(\r))+\r\l^\prime(\r)P(\r)-\r
P^\prime(\r)(2\mu+\l(\r))}{(2\mu+\l(\r))^2}
dx|\\
\di \leq C\int|F|
\left(\f{|F|}{2\mu+\l(\r)}+\f{P(\r)}{2\mu+\l(\r)}\right)P(\r)dx
\di \leq C\left(1+\int \f{|F|^3}{2\mu+\l(\r)}dx\right).
\end{array}
\end{equation}
On the other hand, it holds that
\begin{equation}\label{ue5}
|-\int 2F(u_{1x_2}u_{2x_1}-u_{1x_1}u_{2x_2})dx|\leq C\int|F||\nabla
u|^2 dx.
\end{equation}
Substituting \eqref{ue2}-\eqref{ue5} into \eqref{ue1} yields that
\begin{equation}\label{ue6}
\f12\f{d}{dt}Z^2(t)+\varphi(t)^2\leq
\a\varphi(t)^2+C_\a(Z(t)^2+1)^{2+\f{\v}{1-3\v}}\v^{\f{2}{1-\b}}+C\left[1+\int
\f{|F|^3}{2\mu+\l(\r)}dx+\int|F||\nabla u|^2 dx\right].
\end{equation}
Now it remains to estimate the terms $\di \int
\f{|F|^3}{2\mu+\l(\r)}dx$ and $\di\int|F||\nabla u|^2 dx $ on the
right hand side of \eqref{ue6}. By Lemma \ref{lemma3}, for
$\v\in[0,\f12]$ and $\eta=\v$, it holds that
\begin{equation}\label{Fe}
\|F\|_{2m}\leq C\Big[\|F\|_1+m^{\f12}\|\nabla
F\|_{\f{2m}{m+\v}}^{1-s}\|F\|^s_{2(1-\v)}\Big],
\end{equation}
where $\di s=\f{(1-\v)^2}{m-\v(1-\v)}$ and the positive constant $C$
is independent of $m$ and $\v.$

Choose the positive constant $\v=2^{-m}$ with $m>2$ being integer in
the inequalities \eqref{ue6} and \eqref{Fe}. By the density
estimate \eqref{density-e} in Lemma \ref{lemma-rho}, one can get
\begin{equation}\label{Fe2}
\begin{array}{ll}
\|F\|_1&\di =\int (2\mu+\l(\r))^{-\f12}|F|(2\mu+\l(\r))^{\f12}dx\\
&\di
\leq\left(\f{|F|^2}{2\mu+\l(\r)}\right)^\f12\left(\int(2\mu+\l(\r))dx\right)^{\f12} \leq CZ(t),
\end{array}
\end{equation}
and
\begin{equation}\label{Fe3}
\begin{array}{ll}
\|F\|_{2(1-\v)}^s&\di =\left(\int (2\mu+\l(\r))^{-(1-\v)}|F|^{2(1-\v)}(2\mu+\l(\r))^{1-\v}dx\right)^{\f{s}{2(1-\v)}}\\
&\di
\leq\left(\f{|F|^2}{2\mu+\l(\r)}\right)^{\f s2}\left(\int(2\mu+\l(\r))^{\f{1-\v}{\v}}dx\right)^{\f{s\v}{2(1-\v)}}\\
&\di \leq
CZ(t)^s\left(\|\r\|^{\f{s\b}{2}}_{\f{\b(1-\v)}{\v}}+1\right) \leq
CZ(t)^s\left[\left(\f{\b(1-\v)}{\v}\right)^{\f{s\b}{\b-1}}+1\right]\\&\di
\leq CZ(t)^s\left(\v^{-\f{s\b}{\b-1}}+1\right)
 =
CZ(t)^s\left(2^{\f{ms\b}{\b-1}}+1\right)\leq CZ(t)^s,
\end{array}
\end{equation}
where in the last inequality one has used the fact that
$ms=\f{m(1-\v)^2}{m-\v(1-\v)}\rightarrow1$ as $m\rightarrow +\i.$

Substituting \eqref{fact1} with $r=\f{\v}{m+\v}$, \eqref{Fe2} and
\eqref{Fe3} into \eqref{Fe} yields that
\begin{equation}\label{F-2m}
\begin{array}{ll}
\di \|F\|_{2m}&\di \leq C\Big[Z(t)+m^{\f12}\|\nabla
F\|_{\f{2m}{m+\v}}^{1-s}Z(t)^s\Big]  \leq C\Big[Z(t)+m^{\f12}(\f{m+\v}{\v})^{\f{1-s}{\b-1}}\varphi(t)^{1-s}Z(t)^s\Big]\\
&\di  \leq
C\Big[Z(t)+m^{\f12}(\f{m}{\v})^{\f{1-s}{\b-1}}\varphi(t)^{1-s}Z(t)^s\Big].
\end{array}
\end{equation}
Thus it follows that
\begin{equation}\label{F3}
\begin{array}{ll}
\di \int \f{|F|^3}{2\mu+\l(\r)}dx\di
=\int\f{|F|^{2-\f{1}{m-1}}}{(2\mu+\l(\r))^{1-\f{1}{2(m-1)}}}(\f{1}{2\mu+\l(\r)})^{\f{1}{2(m-1)}}|F|^{1+\f{1}{m-1}}
dx\\
\qquad\quad\di \leq
\int\left(\f{|F|^2}{2\mu+\l(\r)}\right)^{1-\f{1}{2(m-1)}}|F|^{\f{m}{m-1}}
dx\\
\qquad\quad\di \leq \left(\int
\f{|F|^2}{2\mu+\l(\r)}dx\right)^{\f{2m-3}{2(m-1)}}\left(\int
|F|^{2m}dx\right)^{\f{1}{2(m-1)}}\\
\qquad\quad\di \leq Z(t)^{\f{2m-3}{m-1}}\|F\|_{2m}^{\f{m}{m-1}}\leq CZ(t)^{\f{2m-3}{m-1}}\Big[Z(t)+m^{\f12}(\f{m}{\v})^{\f{1-s}{\b-1}}\varphi(t)^{1-s}Z(t)^s\Big]^{\f{m}{m-1}}\\
\qquad\quad\di \leq
C\Big[Z(t)^3+m^{\f{m}{2(m-1)}}(\f{m}{\v})^{\f{(1-s)m}{(\b-1)(m-1)}}\varphi(t)^{\f{(1-s)m}{m-1}}Z(t)^{\f{(2+s)m-3}{m-1}}\Big]\\
\qquad\quad\di \leq
\a\varphi(t)^2+C_\a\Big[Z(t)^3+m^{\f{m}{m(1+s)-2}}(\f{m}{\v})^{\f{2(1-s)m}{(\b-1)(m(1+s)-2)}}Z(t)^{\f{2((2+s)m-3)}{m(1+s)-2}}\Big]\\
\qquad\quad\di \leq
\a\varphi(t)^2+C_\a\Big[(1+Z(t)^2)^2+m(\f{m}{\v})^{\f{2}{\b-1}}(1+Z(t)^2)^{2+\f{1-ms}{m(1+s)-2}}\Big]
\end{array}
\end{equation}
where in the last inequality one has used the fact that
$ms=\f{m(1-\v)^2}{m-\v(1-\v)}\rightarrow1$ with $\v=2^{-m}$ as
$m\rightarrow +\i$.

Furthermore, it holds that
\begin{equation}\label{F-nabla-u}
\begin{array}{ll}
\di \int|F||\nabla u|^2 dx&\di \leq \|F\|_{2m}\|\nabla
u\|_{\f{4m}{2m-1}}^2 \leq C\|F\|_{2m}\Big(\|{\rm div}
u\|_{\f{4m}{2m-1}}^2+\|\o\|_{\f{4m}{2m-1}}^2\Big)\\
&\di \leq
C\|F\|_{2m}\Big(\|\f{F}{2\mu+\l(\r)}\|_{\f{4m}{2m-1}}^2+\|\o\|_{\f{4m}{2m-1}}^2+1\Big).
\end{array}
\end{equation}
Note that
\begin{equation}\label{F-nabla-u1}
\begin{array}{ll}
\di
\|\f{F}{2\mu+\l(\r)}\|_{\f{4m}{2m-1}}^2&\di =\left(\int\f{|F|^{\f{4m}{2m-1}}}{(2\mu+\l(\r))^{\f{4m}{2m-1}}}dx\right)^{\f{2m-1}{2m}}\\
&\di=\left(\int\f{|F|^{\f{2m(2m-3)}{(2m-1)(m-1)}}}{(2\mu+\l(\r))^{\f{4m}{2m-1}}}|F|^{\f{2m}{(2m-1)(m-1)}}dx\right)^{\f{2m-1}{2m}}\\
&\di \leq
\|F\|_{2m}^{\f{1}{m-1}}\left(\int\f{|F|^2}{(2\mu+\l(\r))^{\f{4(m-1)}{2m-3}}}dx\right)^{\f{2m-3}{2(m-1)}}\\
&\di \leq
C\|F\|_{2m}^{\f{1}{m-1}}\left(\int\f{|F|^2}{2\mu+\l(\r)}dx\right)^{\f{2m-3}{2(m-1)}} \leq C\|F\|_{2m}^{\f{1}{m-1}}Z(t)^{\f{2m-3}{m-1}},
\end{array}
\end{equation}
and from $\di\int\o dx=0$, Lemma \ref{lemma2} and \eqref{fact1}, one
has
\begin{equation}\label{F-nabla-u2}
\begin{array}{ll}
\di
\|\o\|_{\f{4m}{2m-1}}^2&\di \leq C\|\o\|_2^{2-\f{1-\v}{m(1-2\v)}}\|\nabla\o\|_{2(1-\v)}^{\f{1-\v}{m(1-2\v)}}\leq C Z(t)^{2-\f{1-\v}{m(1-2\v)}}\Big[\v^{\f{1}{1-\b}}\varphi(t)\Big]^{\f{1-\v}{m(1-2\v)}}\\
&\di \leq C 2^{\f{m(1-\v)}{(\b-1)m(1-2\v)}}Z(t)^{2-\f{1-\v}{m(1-2\v)}}\varphi(t)^{\f{1-\v}{m(1-2\v)}}
\leq C
Z(t)^{2-\f{1-\v}{m(1-2\v)}}\varphi(t)^{\f{1-\v}{m(1-2\v)}}.
\end{array}
\end{equation}
Now substituting \eqref{F-2m}, \eqref{F-nabla-u1} and
\eqref{F-nabla-u2} into \eqref{F-nabla-u} gives that
\begin{equation}\label{F-nabla-u3}
\begin{array}{ll}
&\di \int |F||\nabla u|^2dx\di \leq
C\Big[Z(t)+m^{\f12}(\f{m}{\v})^{\f{1-s}{\b-1}}\varphi(t)^{1-s}Z(t)^s\Big]\Big[1+Z(t)^{2-\f{1-\v}{m(1-2\v)}}\varphi(t)^{\f{1-\v}{m(1-2\v)}}\Big]\\
&\di~~~~~~~~~~~~~~~~~~~~+C\Big[Z(t)+m^{\f12}(\f{m}{\v})^{\f{1-s}{\b-1}}\varphi(t)^{1-s}Z(t)^s\Big]^{1+\f{1}{m-1}}Z(t)^{\f{2m-3}{m-1}}\\
&\di\leq
C\Big[Z(t)+Z(t)^3+Z(t)^{3-\f{1-\v}{m(1-2\v)}}\varphi(t)^{\f{1-\v}{m(1-2\v)}}+m^{\f12}(\f{m}{\v})^{\f{1-s}{\b-1}}\varphi(t)^{1-s}Z(t)^s\\
&\di~+m^{\f12}(\f{m}{\v})^{\f{1-s}{\b-1}}\varphi(t)^{1-s+\f{1-\v}{m(1-2\v)}}Z(t)^{2+s-\f{1-\v}{m(1-2\v)}}+m^{\f12}(\f{m}{\v})^{\f{(1-s)m}{(\b-1)(m-1)}}\varphi(t)^{\f{(1-s)m}{m-1}}Z(t)^{\f{ms+2m-3}{m-1}}\Big]\\
&\di \leq
\a\varphi(t)^2+C_\a\Big[(1+Z^2(t))^2+\Big(m^{\f12}(\f{m}{\v})^{\f{1-s}{\b-1}}Z(t)^s\Big)^{\f{2}{1+s}}\\
&\di~+\Big(m^{\f12}(\f{m}{\v})^{\f{1-s}{\b-1}}Z(t)^{2+s-\f{1-\v}{m(1-2\v)}}\Big)^{\f{2}{1+s-\f{1-\v}{m(1-2\v)}}}+\Big(m^{\f12}(\f{m}{\v})^{\f{(1-s)m}{(\b-1)(m-1)}}Z(t)^{2+\f{ms-1}{m-1}}\Big)^{\f{2(m-1)}{m(s+1)-2}}\Big]\\
&\di \leq
\a\varphi(t)^2+C_\a\Big[(1+Z^2(t))^2+m(\f{m}{\v})^{\f{2}{\b-1}}(1+Z(t)^2)\\
&\di~+m(\f{m}{\v})^{\f{2}{\b-1}}(1+Z(t)^2)^{2+\f{1-ms+(2ms-1)\v}{(1+s)m(1-2\v)-1+\v}}+m(\f{m}{\v})^{\f{2}{\b-1}}(1+Z(t)^2)^{2+\f{1-ms}{m(s+1)-2}}\Big].
\end{array}
\end{equation}
Substituting \eqref{F3} and \eqref{F-nabla-u3} into \eqref{ue6} and
choosing $\a$ sufficiently small yield that
\begin{equation}\label{Z0}
\begin{array}{ll}
\di \f12\f{d}{dt}(Z^2(t))+\f12\varphi(t)^2\leq
C(Z(t)^2+1)^{2+\f{\v}{1-3\v}}\v^{\f{2}{1-\b}}
+ C\Big[(1+Z^2(t))^2+m(\f{m}{\v})^{\f{2}{\b-1}}(1+Z(t)^2)\\
\di~~~~~~~~~~~~~~~~~~~
+m(\f{m}{\v})^{\f{2}{\b-1}}(1+Z(t)^2)^{2+\f{1-ms+(2ms-1)\v}{(1+s)m(1-2\v)-1+\v}}+m^{\f12}(\f{m}{\v})^{\f{2}{\b-1}}(1+Z(t)^2)^{2+\f{1-ms}{m(s+1)-2}}\Big].
\end{array}
\end{equation}
Note that $\lim_{m\rightarrow+\i}[2^m(1-ms)]=2$, and so $1-ms\sim
2\v$ as $m\rightarrow+\i$. Thus for $m$ sufficiently large, one has
$$
\f{1-ms}{m(1+s)-2}\sim \f{2\v}{1-2\v+m-2}=\f{2\v}{m-1-2\v}\leq 4\v,
$$
and
$$
\f{1-ms+(2ms-1)\v}{(1+s)m(1-2\v)-1+\v}=\f{(1-ms)(1-2\v)+\v}{(1+s)m(1-2\v)-1+\v}\sim
\f{3\v-4\v^2}{(m+1-2\v)(1-2\v)-1+\v}\leq 4\v.
$$
Then \eqref{Z0} yields the following inequality for suitably large
$m$,
\begin{equation}\label{Z}
\f12\f{d}{dt}(Z^2(t))+\f12\varphi(t)^2\leq
Cm(\f{m}{\v})^{\f{2}{\b-1}}(1+Z(t)^2)^{2+4\v}.
\end{equation}
Note that
\begin{equation}\label{Z1}
\begin{array}{ll}
Z^2(t)&\di =\int(\mu\o^2+\f{F^2}{2\mu+\l(\r)}) dx\\
&\di \leq C\int[\mu\o^2+(2\mu+\l(\r))({\rm
div}u)^2+\f{P^2(\r)}{2\mu+\l(\r)})] dx\\
&\di\leq C\big(\phi(t)+\int P^2(\r)dx\big)\in L^1(0,T).
\end{array}
\end{equation}
Applying the Gronwall's inequality to \eqref{Z} and using \eqref{Z1} show that
\begin{equation}
\f{1}{(1+Z^2(t))^{4\v}}-\f{1}{(1+Z^2(0))^{4\v}}+Cm\v(\f{m}{\v})^{\f{2}{\b-1}}\geq0.
\end{equation}
Then we have the inequality
\begin{equation}\label{i1}
\f{1}{(1+Z^2(t))^{4\v}}\geq\f{1}{2(1+Z^2(0))^{4\v}},
\end{equation}
provided that
\begin{equation}\label{condition}
Cm\v(\f{m}{\v})^{\f{2}{\b-1}}\leq \f{1}{2(1+Z^2(0))^{4\v}}.
\end{equation}
This condition, \eqref{condition}, is satisfied if
\begin{equation}\label{condition-1}
Cm^{1+\f{2}{\b-1}}2^{-m(1-\f{2}{\b-1})}\leq \f12,
\end{equation}
since
$$
\begin{array}{ll}
\di Z^2(0)&\di=\int\Big[\mu(\o^\d_0)^2+\f{(F^\d_0)^2}{2\mu+\l(\r^\d_0)}\Big] dx\\
&\di \leq C\Big[\|u_0^\d\|^2_{H^2(\mathbb{T}^2)}+\|\r_0^\d\|_{H^3(\mathbb{T}^2)}^{\b}\|u_0^\d\|^2_{H^2(\mathbb{T}^2)}+\|\r_0^\d\|_{H^3(\mathbb{T}^2)}^{2\g}\Big]\leq C.
\end{array}
$$
Now if $\b>3,$ that is, $1-\f{2}{\b-1}>0$, then we can choose
sufficiently large $m>2$ to guarantee the condition
\eqref{condition-1}. Consequently, the inequality \eqref{i1} is
satisfied with  $\b>3$ and sufficiently large $m>2$. Then
\begin{equation}\label{Z2}
Z^2(t)\leq 2^{2^{m-1}}(1+Z^2(0))-1\leq C,
\end{equation}
and \begin{equation}\label{Z3}
 \int_0^T\varphi(t)dt\leq C.
\end{equation}
Thus the proof of Lemma \ref{lemma-u-der} is completed. $\hfill\Box$

\underline{Step 5: Second order derivative estimates for the
velocity:}

\begin{Lemma}\label{lemma-u-sec}
  There exists a positive constant $C$ independent of $\d$, such that
  \begin{equation}
    \sup_{t\in[0,T]}\int\r(H^2+L^2)
     dx+\int_0^T\int\mu(H_{x_1}-L_{x_2})^2+(2\mu+\l(\r))(H_{x_2}+L_{x_1})^2dxdt\leq C.
  \end{equation}
\end{Lemma}
{\bf Proof:} Multiplying the equations, $\eqref{H-L}_1$ and $\eqref{H-L}_2$, by $H$ and $L$, respectively, summing the resulted equations together, and integrating with respect to $x$ over $\mathbb{T}^2$ lead to
\begin{equation}\label{HL1}
\begin{array}{ll}
\di
\f12\f{d}{dt}\int\r(H^2+L^2)dx+\int\mu(H_{x_1}-L_{x_2})^2+(2\mu+\l(\r))(H_{x_2}+L_{x_1})^2dx
dx\\
\di =\int\r(H^2+L^2){\rm div}u dx-\int\mu\o{\rm div}u (L_{x_2}-H_{x_1}) dx\\
\di~~ -\int\r(2\mu+\l(\r))\big[F(\f{1}{2\mu+\l(\r)})^\prime+(\f{P(\r)}{2\mu+\l(\r)})^\prime\big]{\rm div}u(H_{x_2}+L_{x_1}) dx\\
\di~~-\int \big[H(u_{x_2}\cdot\nabla F+\mu
u_{x_1}\cdot\nabla\o)+L(u_{x_1}\cdot\nabla F-\mu u_{x_2}\cdot\nabla\o)\big]dx\\
\di~~+\int(2\mu+\l(\r))[(u_{1x_1})^2+2u_{1x_2}u_{2x_1}+(u_{2x_2})^2](H_{x_2}+L_{x_1})
dx.
\end{array}
\end{equation}
Set
\begin{equation}
Y(t)=\left(\int\r(H^2+L^2)dx\right)^{\f12},
\end{equation}
and
\begin{equation}
\psi(t)=\left(\int\mu(H_{x_1}-L_{x_2})^2+(2\mu+\l(\r))(H_{x_2}+L_{x_1})^2dx\right)^{\f12}.
\end{equation}
Note that
$$
\begin{array}{ll}
\di \int(|\nabla H|^2+|\nabla L|^2)dx&\di =\int
(H_{x_1}^2+H_{x_2}^2+L_{x_1}^2+L_{x_2}^2)dx\\
&\di =\int\big[(H_{x_1}-L_{x_2})^2+(H_{x_2}+L_{x_1})^2\big] dx\leq \f{1}{\mu}\psi^2(t).
\end{array}
$$
Thus it holds that
\begin{equation}\label{fact3}
\|\nabla(H,L)\|_2(t)\leq C\psi(t),\qquad\forall t\in[0,T].
\end{equation}
Then it follows from the elliptic system
$$
\mu\o_{x_1}+F_{x_2}=\r H,\qquad\qquad -\mu\o_{x_2}+F_{x_1}=\r L,
$$
 that
\begin{equation}\label{fact4}
\|\nabla(F,\o)\|_p\leq C\|\r(H,L)\|_p, \qquad\forall 1<p<+\i.
\end{equation}
Furthermore, since $\di\int (\mu\o_{x_1}+F_{x_2}) dx=0,$ by the mean
value theorem, there exists a point $x_*\in\mathbb{T}^2$, such that
$(\mu\o_{x_1}+F_{x_2})(x_*,t)=0,$ and so $H(x_*,t)=0.$ Similarly,
there exists a point $x_*^\prime$, such that $L(x_*^\prime,t)=0.$
Therefore, by the Poincare inequality, it holds that
\begin{equation}\label{fact5}
\|(H,L)\|_{p}\leq C\|\nabla(H,L)\|_2, \qquad \forall 1\leq p<+\i,
\end{equation}
where $C$ may depend on $p$.

 Now we estimate the right hand side of
\eqref{HL1} term by term. First, by the H${\rm\ddot{o}} $lder
inequality, \eqref{fact5} and the density estimate
\eqref{density-e},  it holds that
\begin{equation}\label{HL2}
\begin{array}{ll}
\di |\int\r(H^2+L^2){\rm div} u dx| =|\int\r(H^2+L^2) \f{F+P(\r)}{2\mu+\l(\r)}dx|\\
\di\qquad\quad \leq\|\sqrt\r(H,L)\|_2\|(H,L)\|_4\|\f{\sqrt\r(F+P(\r))}{2\mu+\l(\r)}\|_4 \leq CY(t)\psi(t)(1+\|F\|_4).
\end{array}
\end{equation}
Note that
\begin{equation}\label{F4}
\begin{array}{ll}
\|(F,\o)\|_4 \leq C(\|\nabla (F,\o)\|_{\f32}+\|(F,\o)\|_1)\\
\di~\leq C\Big[\|\nabla
(F,\o)\|_{\f32}+\Big(\int\f{F^2}{2\mu+\l(\r)}dx\Big)^{\f12}\Big(\int(2\mu+\l(\r))dx\Big)^{\f12}+\|\o\|_2\Big] \leq C\Big[Y(t)+1\Big],
\end{array}
\end{equation}
where in the last inequality one has used the estimate \eqref{fact1}
with $r=\f14$ and the estimate \eqref{Z2}.

Substituting \eqref{F4} into \eqref{HL2} yields that
\begin{equation}\label{HL21}
\begin{array}{ll}
\di |\int\r(H^2+L^2){\rm div} u dx|&\di \leq CY(t)\psi(t)(Y(t)+1) \leq \a \psi^2(t)+C_\a (Y(t)+1)^4.
\end{array}
\end{equation}
Second, direct estimates give
\begin{equation}\label{HL22}
\begin{array}{ll}
\di |-\int\mu\o{\rm div}u (L_{x_2}-H_{x_1}) dx| \leq \mu\left(\int
(L_{x_2}-H_{x_1})^2dx\right)^{\f12}\left(\int\o^2({\rm div}
u)^2dx\right)^{\f12}\\[2mm]
\di \leq \a\psi^2(t)+C_\a \int\o^2({\rm div} u)^2dx\leq \a\psi^2(t)+C_\a \|\o\|_4^2\|\f{F+P(\r)}{2\mu+\l(\r)}\|_4^2\\
\di \leq \a\psi^2(t)+C_\a \|\o\|_4^2(1+\|F\|_4^2) \leq \a \psi^2(t)+C_\a (Y(t)+1)^4.
\end{array}
\end{equation}
Similarly, one has
\begin{equation}\label{HL23}
\begin{array}{ll}
\di
|-\int\r(2\mu+\l(\r))\big[F(\f{1}{2\mu+\l(\r)})^\prime+(\f{P(\r)}{2\mu+\l(\r)})^\prime\big]{\rm
div}u(H_{x_2}+L_{x_1}) dx|\\
\di \leq \a\int
(2\mu+\l(\r))(H_{x_2}+L_{x_1})^2dx\\
\di \qquad\qquad+C_\a\int\r^2(2\mu+\l(\r))\big[F(\f{1}{2\mu+\l(\r)})^\prime+(\f{P(\r)}{2\mu+\l(\r)})^\prime\big]^2({\rm
div}u)^2 dx\\
\di \leq
\a\psi^2(t)+C_\a\int\r^2\big[F(\f{1}{2\mu+\l(\r)})^\prime+(\f{P(\r)}{2\mu+\l(\r)})^\prime\big]^2\f{|F|^2+P^2(\r)}{2\mu+\l(\r)} dx\\
\di \leq \a\psi^2(t)+C_\a (1+\|F\|_4^4)\leq \a \psi^2(t)+C_\a (Y(t)+1)^4.
\end{array}
\end{equation}
Next,
\begin{equation}\label{HL24}
\begin{array}{ll}
\di |-\int \big[H(u_{x_2}\cdot\nabla F+\mu
u_{x_1}\cdot\nabla\o)+L(u_{x_1}\cdot\nabla F-\mu u_{x_2}\cdot\nabla\o)\big]dx|\\
\di\leq C\int|(H,L)||\nabla u||\nabla(F,\o)| dx\\
\di\leq C\|(H,L)\|_8\|\nabla u\|_2\|\nabla(F,\o)\|_{\f83}\leq C\|\nabla(H,L)\|_2\|\r(H,L)\|_{\f83},
\end{array}
\end{equation}
where one has used the fact that
$$
\|\nabla u\|_2\leq C(\|{\rm div}u\|_2+\|\o\|_2)\leq
C(\|\f{F+P(\r)}{2\mu+\l(\r)}\|_2+\|\o\|_2)\leq C.
$$
Note that
\begin{equation}\label{HL25}
\begin{array}{ll}
\di \|\r(H,L)\|_{\f83}&\di=\left(\int\r^{\f83}|(H,L)|^{\f83}dx\right)^{\f38} =\left(\int\sqrt\r|(H,L)||(H,L)|^{\f53}\r^{\f{13}{6}}dx\right)^{\f38}\\
&\di\leq
\|\sqrt\r(H,L)\|_2^{\f38}\|(H,L)\|_4^{\f58}\|\r\|_{26}^{\f{13}{16}}
\leq CY(t)^{\f38}\|\nabla(H,L)\|_2^{\f58}.
\end{array}
\end{equation}
It follows from \eqref{HL24} and \eqref{HL25} that
\begin{equation}\label{HL26}
\begin{array}{ll}
\di |-\int \big[H(u_{x_2}\cdot\nabla F+\mu
u_{x_1}\cdot\nabla\o)+L(u_{x_1}\cdot\nabla F-\mu u_{x_2}\cdot\nabla\o)\big]dx|\\
\di\leq CY(t)^{\f38}\|\nabla(H,L)\|_2^{\f{13}8}\leq
CY(t)^{\f38}\psi(t)^{\f{13}8} \leq \a \psi(t)^2+C_\a Y(t)^2.
\end{array}
\end{equation}
Moreover,
\begin{equation}\label{HL27}
\begin{array}{ll}
\di|\int(2\mu+\l(\r))[(u_{1x_1})^2+2u_{1x_2}u_{2x_1}+(u_{2x_2})^2](H_{x_2}+L_{x_1})
dx|\\
\di\leq
\a\psi(t)^2+C_\a\int(2\mu+\l(\r))[(u_{1x_1})^2+2u_{1x_2}u_{2x_1}+(u_{2x_2})^2]^2dx\\
\di\leq
\a\psi(t)^2+C_\a\|2\mu+\l(\r)\|_2\|\nabla u\|_8^4\\
\di\leq
\a\psi(t)^2+C_\a(\|{\rm div} u\|_8^4+\|\o\|_8^4)\leq
\a\psi(t)^2+C_\a(\|(F,\o)\|_8^4+1)\\
\di\leq
\a\psi(t)^2+C_\a(\|\nabla(F,\o)\|_{\f85}^4+1) \leq \a \psi(t)^2+C_\a (1+Y(t))^4.
\end{array}
\end{equation}
Substituting the estimates \eqref{HL21}-\eqref{HL23}, \eqref{HL26}
and \eqref{HL27} into \eqref{HL1}, one can arrive at
\begin{equation}
\f12\f{d}{dt}(Y^2(t))+\psi^2(t)\leq 5\a\psi^2(t)+C_\a(1+Y^2(t))^2.
\end{equation}
Choosing $5\a=\f12$, noting that $Y^2(t)=\varphi^2(t)\in
L^1(0,T)$, and then using Gronwall's inequality yield that
\begin{equation}\label{Y-e}
Y^2(t)+\int_0^T\psi^2(t)dt\leq Y^2(0)+C.
\end{equation}
Now we calculate the initial values $Y^2(0).$ By the approximate compatibility
condition \eqref{new-cc}, one has
$$
\mathcal{L}_{\r^\d_0}u^\d_0-\nabla P^\d_0=\sqrt\r_0 g, ~~{\rm with}~~ g\in L^2(\mathbb{T}^2).
$$
On the other hand, it holds that
\begin{equation}
\begin{array}{ll}
\mathcal{L}_{\r^\d_0}u^\d_0&\di =\mu\Delta u^\d_0+\nabla((\mu+\l(\r^\d_0)){\rm
div}u^\d_0) =\mu\Delta u^\d_0+\nabla(F^\d_0-\mu{\rm div}u^\d_0+P^\d_0)\\
&\di=[\mu\nabla({\rm div}u^\d_0)-\mu\nabla\times(\nabla\times
u^\d_0)]+\nabla(F^\d_0-\mu{\rm div}u^\d_0+P^\d_0)
\end{array}
\end{equation}
where $F_0^\d=(2\mu+\l(\r_0^\d)){\rm div} u_0^\d-P_0^\d$ and similarly one can define
 $\o_0^\d, L_0^\d, H_0^\d$,  $\nabla\times$ denotes the  3-dimensional {\it curl} operator,
and
$$
\nabla\times(\nabla\times
u^\d_0)=(\partial_{x_2}\o^\d_0,-\partial_{x_1}\o^\d_0,0)
$$
is regarded as the 2-dimensional vector
$(\partial_{x_2}\o^\d_0,-\partial_{x_1}\o^\d_0)^t$.

Thus
\begin{equation}\label{82}
\begin{array}{ll}
\mathcal{L}_{\r^\d_0}u^\d_0-\nabla P^\d_0&\di =\nabla
F^\d_0-\mu(\partial_{x_2}\o^\d_0,-\partial_{x_1}\o^\d_0)^t\\
&\di
=(F^\d_{0x_1}-\mu\partial_{x_2}\o^\d_0,F^\d_{0x_2}+\mu\partial_{x_1}\o^\d_0)^t=\r^\d_0(L^\d_0,H^\d_0)^t.
\end{array}
\end{equation}
Therefore
\begin{equation}
\sqrt{\r_0}g=\r_0^\d(L^\d_0,H^\d_0)^t.
\end{equation}
Consequently, it holds that
\begin{equation}\label{id1}
Y^2(0)=\|\sqrt{\r^\d_0}(L^\d_0,H^\d_0)\|_2^2=\|\f{\sqrt{\r_0}}{\sqrt{\r_0^\d}} g\|_2^2\leq C.
\end{equation}
This, together with \eqref{Y-e}, shows that
\begin{equation}
Y^2(t)+\int_0^T \psi^2(t)dt\leq C.
\end{equation}
This completes the proof of Lemma \ref{lemma-u-sec}. $\hfill\Box$

\begin{Remark}\label{remark31}
Similar to the derivation of \eqref{82}, one can get that for any $t\in[0,T]$,
$$
\mathcal{L}_\r u-\nabla P(\r)=\r(L,H)^t.
$$
Then it follows from the momentum equation $\eqref{CNS}_2$ that
\begin{equation}\label{ut-HL}
u_t=(L,H)^t-u\cdot \nabla u.
\end{equation}
The above identity can also be obtained directly from \eqref{ns1}.
\end{Remark}

\underline{Step 6. Upper bound of the density:} We are now ready to derive the upper
bound for the density in the super-norm independent of $\d$, which is crucial for the
proof of Theorem \ref{theorem2} as in \cite{hx2, hlx1, hlx}. First, we have

\begin{Lemma}\label{lemma-F-o}
  It holds that
  \begin{equation}
\int_0^T\|(F,\o)\|_\i^3dt\leq C.
  \end{equation}
\end{Lemma}
{\bf Proof:} By \eqref{fact4} with $p=3$, one has
 \begin{equation}
\begin{array}{ll}
\di \int_0^T\|\nabla(F,\o)\|_3^3dt&\di \leq C\int_0^T\|\r(H,L)\|_3^3dt =C\int\int\r^3|(H,L)|^3dxdt\\
&\di =C\int\int\sqrt\r|(H,L)||(H,L)|^2\r^{\f52}dxdt\\
&\di \leq C\int\|\sqrt\r(H,L)\|_2\|(H,L)\|_8^2\|\r\|^{\f52}_{10}dt\\
&\di\leq C\int_0^T\|\nabla(H,L)\|_2^2dt\leq C\int_0^T \psi^2(t)\leq
C,
\end{array}
\end{equation}
which, combined with the estimates in Lemma \ref{lemma3}, yields that
\begin{equation}\label{F-omega-infty}
\int_0^T\|(F,\o)\|_\i^3dt\leq
\int_0^T\|(F,\o)\|_{W^{1,3}(\mathbb{T}^2)}^3dt\leq C.
\end{equation}
The proof of Lemma \ref{lemma-F-o} is finished. $\hfill\Box$

With Lemma \ref{lemma-F-o} in hand, we can obtain the uniform upper bound for the density.

\begin{Lemma}
  It holds that
  \begin{equation}
\r(t,x)\leq C,\qquad\forall (t,x)\in [0,T]\times \mathbb{T}^2.
  \end{equation}
\end{Lemma}
{\bf Proof:} From the continuity equation $\eqref{CNS}_1$, we have
\begin{equation}
\t(\r)_t+u\cdot\nabla\t(\r)+P(\r)+F=0,
\end{equation}
where $\t(\r)$ is defined in \eqref{theta}.

Along the particle path $\vec{X}(\tau;t,x)$ through the point
$(t,x)\in[0,T]\times\mathbb{T}^2$ defined by
\begin{equation}
\left\{
\begin{array}{ll}
\di \f{d\vec{X}(\tau;t,x)}{d\tau}=u(\tau,\vec{X}(\tau;t,x)),\\
 \di
\vec{X}(\tau;t,x)|_{\tau=t}=x,
\end{array}
\right.
\end{equation}
there holds the following ODE
\begin{equation}
\f{d}{d\tau}\t(\r)(\tau,\vec{X}(\tau;t,x))=-P(\r)(\tau,\vec{X}(\tau;t,x))-F(\tau,\vec{X}(\tau;t,x)),
\end{equation}
which is integrated over $[0,t]$ to yield that
\begin{equation}\label{theta-1}
\t(\r)(t,x)-\t(\r_0)(\vec{X}_0)=-\int_0^t(P(\r)+F)(\tau,\vec{X}(\tau;t,x))d\tau,
\end{equation}
with $\vec{X}_0=\vec{X}(\tau;t,x)|_{\tau=0}$.

It follows from \eqref{theta-1} that
\begin{equation}
2\mu\ln\f{\r(t,x)}{\r_0(\vec{X}_0)}+\f{1}{\b}\r^\b(t,x)+\int_0^tP(\r)(\tau,\vec{X}(\tau;t,x))d\tau=\f1\b\r_0(\vec{X}_0)^{\b}-\int_0^tF(\tau,\vec{X}(\tau;t,x))d\tau.
\end{equation}
So
\begin{equation}
2\mu\ln\f{\r(t,x)}{\r_0(\vec{X}_0)}\leq
\f1\b\|\r_0\|_\i^{\b}+\int_0^t\|F(\tau,\cdot)\|_\i d\tau\leq C,
\end{equation}
which implies that
$$
\f{\r(t,x)}{\r_0(\vec{X}_0)}\leq C.
$$
Therefore, we have
\begin{equation}
\r(t,x)\leq C, \qquad \forall (t,x)\in[0,T]\times\mathbb{T}^2.
\end{equation}
Hence the Lemma is proved. $\hfill\Box$

As an immediate consequence of the upper bound of the density, one has
\begin{Lemma}\label{lemma3.9}
  It holds that for any $1< p<\i$,
  \begin{equation}
    \int_0^T\big(\|{\rm div}u\|_\i^3+\|\nabla(F,\o)\|_p^2\big)dt\leq C.
  \end{equation}
  \end{Lemma}
{\bf Proof:} First, note that
\begin{equation}\label{div-u-infty}
\int_0^T\|{\rm div}u\|_\i^3dt\leq C\int_0^T
(\|F\|_\i^3+\|P(\r)\|_\i^3)dt\leq C.
\end{equation}
Then for any $1< p<\i$,
\begin{equation}\label{nabla-F-omega}
\begin{array}{ll}
\di \int_0^T\|\nabla(F,\o)\|_p^2dt&\di \leq C\int_0^T\|\r(H,L)\|_p^2dt\\
&\di\leq C\int_0^T\|(H,L)\|_p^2dt \leq C\int_0^T\|\nabla(H,L)\|_2^2dt\leq C.
\end{array}
\end{equation}
Thus Lemma \ref{lemma3.9} is proved. $\hfill\Box$

\section{Higher order estimates}
\setcounter{equation}{0}

With the approximate solutions and basic estimates at hand, we can derive
some uniform estimates on their higher order derivatives easily as in \cite{hx2, hlx1,hlx}.
We start with estimates on first order derivatives.

\begin{Lemma}\label{lemma4.1}
  It holds that for any $1\leq p<+\i$,
  \begin{equation}
    \sup_{t\in[0,T]}\|(\nabla\r,\nabla P(\r))(t,\cdot)\|_p+\int_0^T\|\nabla u\|_\i^2dt\leq C.
  \end{equation}
  \end{Lemma}
{\bf Proof:} Applying the operator $\nabla$ to the continuity equation
$\eqref{CNS}_1$, one has
\begin{equation}\label{nabla-rho}
(\nabla\r)_t+\nabla
u\cdot\nabla\r+u\cdot\nabla(\nabla\r)+\nabla\r{\rm div}
u+\r\nabla({\rm div}u)=0.
\end{equation}
Multiplying the equation \eqref{nabla-rho} by
$p|\nabla\r|^{p-2}\nabla\r$ with $p\geq2$ implies that
\begin{equation}
(|\nabla\r|^p)_t+{\rm div}(u|\nabla\r|^p)+(p-1)|\nabla\r|^p{\rm
div}u+p|\nabla\r|^{p-2}\nabla\r\cdot(\nabla
u\cdot\nabla\r)+p\r|\nabla\r|^{p-2}\nabla\r\cdot\nabla({\rm div}
u)=0.
\end{equation}
Integrating over $\mathbb{T}^2$ gives
\begin{equation}
\begin{array}{ll}
\di \f{d}{dt}\|\nabla\r\|_p^p\\
\di=-(p-1)\int |\nabla u|^p{\rm div}
udx-p\int|\nabla\r|^{p-2}\nabla\r\cdot (\nabla
u\cdot\nabla\r)dx-p\int\r|\nabla\r|^{p-2}\nabla\r\cdot\nabla({\rm
div} u)dx\\
\di\leq (p-1)\|{\rm div}u\|_\i\|\nabla\r\|_p^p+p\|\nabla
u\|_\i\|\nabla\r\|_p^p+p\|\r\|_\i\|\nabla\r\|_p^{p-1}\|\nabla{\rm
div}u\|_p.
\end{array}
\end{equation}
This implies that
\begin{equation}\label{100}
\begin{array}{ll}
\di \f{d}{dt}\|\nabla\r\|_p&\di \leq C\Big[\|\nabla
u\|_\i\|\nabla\r\|_p+\|\nabla{\rm div}u\|_p\Big]\\
&\di \leq C\Big[\|\nabla
u\|_\i\|\nabla\r\|_p+\|\nabla\big(\f{F+P(\r)}{2\mu+\l(\r)}\big)\|_p\Big]\\
&\di \leq C\Big[\Big(\|\nabla
u\|_\i+\|F\|_\i+1\Big)\|\nabla\r\|_p+\|\nabla F\|_p\Big].
\end{array}
\end{equation}
By Remark \ref{remark31}, one has
\begin{equation}\label{elliptic-HL}
\mathcal{L}_\r u=\nabla P(\r)+\r(L,H)^t.
\end{equation}
Thus the elliptic estimates and \eqref{fact5} yields that for any $1<p<\i,$
\begin{equation}\label{e5}
\begin{array}{ll}
\|\nabla^2u\|_{p}&\di \leq C\big[\|\nabla P(\r)\|_p+\|\r(L,H)\|_p\big]\\
&\di\leq C\big[\|\nabla \r\|_p+\|(L,H)\|_p\big]\leq C\big[\|\nabla \r\|_p+\|\nabla(L,H)\|_2\big].
\end{array}
\end{equation}
By Beal-Kato-Majda type inequality (see \cite{hlx}-\cite{hx2} or \cite{Kazhikhov}), it holds that
\begin{equation}\label{103}
\begin{array}{ll}
\di\|\nabla u\|_\i&\di \leq C\big(\|{\rm
div}u\|_\i+\|\o\|_\i\big)\ln(e+\|\nabla^2u\|_3)\\
&\di\leq C\big(\|{\rm
div}u\|_\i+\|\o\|_\i\big)\ln(e+\|\nabla\r\|_3)+C\big(\|{\rm
div}u\|_\i+\|\o\|_\i\big)\ln(e+\|\nabla(H,L)\|_2).
\end{array}
\end{equation}
The combination of \eqref{100} with $p=3$ and \eqref{103} yields
that
\begin{equation}
\begin{array}{ll}
\di\f{d}{dt}\|\nabla\r\|_3&\di \leq C\Big[\big(\|{\rm
div}u\|_\i+\|\o\|_\i\big)\ln(e+\|\nabla(H,L)\|_2)+\|F\|_\i+1\Big]\|\nabla\r\|_3\\
&\di~~~+C\big(\|{\rm div}u\|_\i+\|\o\|_\i\big)\|\nabla
\r\|_3\ln(e+\|\nabla\r\|_3)+C\|\nabla F\|_3.
\end{array}
\end{equation}
By the estimates \eqref{F-omega-infty}, \eqref{div-u-infty},
\eqref{nabla-F-omega} and the Gronwall's inequality, it holds that
\begin{equation}
\sup_{t\in[0,T]}\|\nabla\r\|_3\leq C,
\end{equation}
which, together with \eqref{F-omega-infty}, \eqref{div-u-infty}, \eqref{e5} and
\eqref{103},  yields that
\begin{equation}\label{nabla-u-infty}
\int_0^T\|\nabla u\|_\i^2dt\leq C.
\end{equation}
Therefore, by \eqref{nabla-u-infty}, Lemma \ref{lemma-F-o}, Lemma
\ref{lemma3.9} and Gronwall inequality, one can derive from
\eqref{100} that
\begin{equation}\label{nabla-rho-p}
\sup_{t\in[0,T]}\|\nabla\r\|_p\leq C(\|\nabla\r_0\|_p+1),
\qquad\forall p\in[1,+\i).
\end{equation}
Thus the proof of Lemma \ref{lemma4.1} is completed. $\hfill\Box$

\begin{Lemma}\label{lemma-r-h2}
  It holds that for any $1\leq p<+\i$,
  \begin{equation}
    \sup_{t\in[0,T]}\Big[\|u(t,\cdot)\|_\i+\|\nabla u\|_p+\|(\r_t,P_t)\|_p+\|(\r_t,P(\r)_t)\|_{H^1}+\|(\r,u)\|_{H^2}\Big]+\int_0^T\|u\|_{H^3}^2dt\leq C.
  \end{equation}
  \end{Lemma}
{\bf Proof:} By $L^2-$estimates to the elliptic system \eqref{elliptic-HL},
one has
\begin{equation}
\begin{array}{ll}
\di\sup_{t\in[0,T]}\|u\|_{H^2}&\di \leq C
\sup_{t\in[0,T]}\big(\|\nabla
P(\r)\|_2+\|\r(H,L)\|_2\big)\\
&\di \leq
C\sup_{t\in[0,T]}\big(\|\nabla P(\r)\|_2+\|\sqrt\r(H,L)\|_2\big)\leq
C.
\end{array}
\end{equation}
It follows from the  Sobolev embedding theorem that
\begin{equation}\label{u-infty}
\sup_{[0,T]\times\mathbb{T}^2}|u(t,x)|\leq C,\qquad
\sup_{t\in[0,T]}\|\nabla u\|_p\leq C, ~~\forall 1\leq p<+\i.
\end{equation}
Due to $\eqref{CNS}_1$, one can get
$\r_t=-u\cdot\nabla\r-\r~{\rm div}u$ and $P_t=-u\cdot\nabla P-\r
P^\prime(\r)~{\rm div}u$, which, together with the uniform upper
bound of the density  and the estimates in Lemma \ref{lemma4.1}
and \eqref{u-infty}, yields that
\begin{equation}\label{rho-t-p}
\sup_{t\in[0,T]}\|(\r_t,P_t)\|_p\leq C, \qquad\forall p\in[1,+\i).
\end{equation}
Applying $\nabla^2$ to the continuity equation $\eqref{CNS}_1$, then
multiplying the resulted equation by $\nabla^2\r$, and then
integrating over the torus $\mathbb{T}^2$, one can get that
\begin{equation}\label{nabla-rho-2}
\begin{array}{ll}
\di\f{d}{dt}\|\nabla^2\r\|_2^2&\di \leq C\Big[\|\nabla u\|_\i\|\nabla^2\r\|^2_2+\|\nabla\r\|_4\|\nabla^2 \r\|_2\|\nabla^2 u\|_4+\|\r\|_\i\|\nabla^2\r\|_2\|\nabla^3u\|_2\Big]\\
&\di\leq C\Big[\big(\|\nabla
u\|_\i+1)\|\nabla^2\r\|^2_2+\|\nabla^3u\|^2_2+1\Big].
\end{array}
\end{equation}
Similarly,
\begin{equation}\label{P2}
\f{d}{dt}\|\nabla^2P(\r)\|_2^2\leq C\Big[\big(\|\nabla
u\|_\i+1)\|\nabla^2P(\r)\|^2_2+\|\nabla^3u\|^2_2+1\Big].
\end{equation}
Note that \eqref{elliptic-HL} implies that
$$
\mathcal{L}_\r (\nabla
u)=\nabla^2P(\r)+\nabla[\r(H,L)]+\nabla(\nabla\l(\r){\rm
div}u):=\Phi.
$$
Then the standard elliptic estimates give that
\begin{equation}\label{nabla-u-3}
\begin{array}{ll}
\di \|u\|_{H^3}&\di \leq C\Big[\|u\|_{H^1}+\|\Phi\|_2\Big]\\
&\di \leq C\Big[\|u\|_{H^1}+\|\nabla^2P(\r)\|_2+\|\r\|_\i\|\nabla(H,L)\|_2+\|\nabla\r\|_4\|(H,L)\|_4\\
&\di ~~~~~+\|\nabla^2\r\|_2\|{\rm
div}u\|_\i+\|\nabla\r\|_4\|\nabla^2u\|_4\Big],
\end{array}
\end{equation}
and
$$
\|\nabla^2u\|_4\le C \|\nabla P (\rho)\|_4+\|\rho(L,H)\|_4\le
C(1+\|\nabla(H,L)\|_2).
$$
Consequently,
\begin{equation}\label{nabla-u-30}
\di \|u\|_{H^3}\di \leq
C\Big[1+\|\nabla^2P(\r)\|_2+\|\nabla(H,L)\|_2+\|\nabla^2\r\|_2\|{\rm
div}u\|_\i\Big].
\end{equation}
Substituting \eqref{nabla-u-3} into \eqref{nabla-rho-2} and
\eqref{P2} yields that
\begin{equation}
\f{d}{dt}\|(\nabla^2\r,\nabla^2P(\r))\|_2^2\leq C\Big[\big(\|\nabla
u\|^2_\i+1\big)\|(\nabla^2\r,\nabla^2P(\r))\|^2_2+\|\nabla(H,L)\|^2_2+1\Big].
\end{equation}
Then the Gronwall's inequality yields that
\begin{equation}
\begin{array}{ll}
\di\|(\nabla^2\r,\nabla^2P(\r))\|_2^2(t)&\di \leq
\Big(\|(\nabla^2\r_0,\nabla^2P_0)\|_2^2+C\int_0^T(\|\nabla(H,L)\|_2^2+1)dt\Big)e^{\di
C\int_0^T\big(\|\nabla u\|^2_\i+1\big)dt}\\
&\di \leq C,
\end{array}
\end{equation}
which also implies that
\begin{equation}
\sup_{t\in[0,T]}\big(\|(\r,P(\r))\|_{H^2}+\|(\r_t,P(\r)_t)\|_{H^1}\big)+
\int_0^T\|u\|_{H^3}^2dt\leq C.
\end{equation}
The proof of Lemma \ref{lemma-r-h2} is completed. $\hfill\Box$

\begin{Lemma}\label{nabla-ut-1}
It holds that
  \begin{equation}
\sup_{t\in[0,T]}\|\sqrt\r u_t\|_2^2(t)+\int_0^T\|u_t\|_{H^1}^2dt\leq C.
\end{equation}
\end{Lemma}
{\bf Proof:} The momentum equation $\eqref{CNS}_2$ can be written as
\begin{equation}\label{moment-1}
\r u_t+\r u\cdot\nabla u+\nabla P(\r)=\mathcal{L}_\r u:=\mu\Delta
u+\nabla((\mu+\l(\r)){\rm div} u).
\end{equation}
Applying $\partial_t$ to the above equation gives that
\begin{equation}\label{u-tt}
\r u_{tt}+\r u\cdot\nabla u_t+\nabla P(\r)_t=\mu\Delta
u_t+\nabla((\mu+\l(\r)){\rm div} u_t)-\r_tu_t-\r_tu\cdot \nabla u-\r
u_t\cdot\nabla u+\nabla(\l(\r)_t{\rm div}u).
\end{equation}
Multiplying the equation \eqref{u-tt} by $u_t$ and integrating the
resulting equation with respect to $x$ over $\mathbb{T}^2$
imply that
\begin{equation}\label{u-tt1}
\begin{array}{ll}
\di\f12\f{d}{dt}\int \r|u_t|^2dx+\int\Big(\mu|\nabla
u_t|^2+(\mu+\l(\r))|{\rm div}
u_t|^2\Big)dx\\
\di=-\int \nabla P(\r)_t\cdot u_tdx-\int\r_t|u_t|^2 dx-\int
\r_t(u\cdot\nabla u)\cdot u_tdx\\
\di ~~-\int\r(u_t\cdot\nabla u)\cdot u_t dx+\int\nabla(\l(\r)_t{\rm
div}u)\cdot u_t dx.
\end{array}
\end{equation}
Notice that
\begin{equation}
\begin{array}{ll}
\di -\int \nabla P(\r)_t\cdot u_tdx&\di=\int P(\r)_t {\rm div }u_t dx\\
&\di\leq \f{\mu}{4}\int|{\rm div }u_t|^2 dx+C\int |P_t|^2dx\leq \f{\mu}{4}\int|{\rm div }u_t|^2 dx+C,
\end{array}
\end{equation}
\begin{equation}
\begin{array}{ll}
\di -\int\r_t|u_t|^2 dx=\int {\rm div }(\r u)|u_t|^2 dx=-2\int\r(u\cdot\nabla u_t)\cdot u_tdx\\
\qquad\qquad\di\leq \f{\mu}{8}\int|\nabla u_t|^2
dx+C\|u\|_\i^2\|\sqrt\r\|_\i^2\|\sqrt\r u_t\|_2^2\leq \f{\mu}{4}\int|\nabla u_t|^2 dx+C\|\sqrt\r u_t\|_2^2,
\end{array}
\end{equation}
\begin{equation}
\begin{array}{ll}
\di -\int \r_t(u\cdot\nabla u)\cdot u_tdx=\int {\rm div }(\r u)[(u\cdot\nabla u)\cdot u_t]dx=-\int\r u\cdot\nabla[(u\cdot\nabla u)\cdot u_t]dx\\
\di\qquad\quad\leq \|\r\|_\i\|u\|_\i^2\|\nabla u_t\|_2\|\nabla u\|_2+\|u\|_\i\|\sqrt\r\|_\i\|\sqrt\r u_t\|_2\big(\|\nabla u\|_4^2+\|u\|_\i\|\nabla^2u\|_2\big)\\
\di\qquad\quad\leq \f{\mu}{4}\int|\nabla u_t|^2 dx+C\big(\|\sqrt\r
u_t\|_2^2+\|(\nabla u,\nabla^2u)\|_2^2+\|\nabla u\|_4^4\big)\\
\di\qquad\quad\leq \f{\mu}{4}\int|\nabla u_t|^2 dx+C\big(\|\sqrt\r
u_t\|_2^2+1\big),
\end{array}
\end{equation}
\begin{equation}
|-\int\r(u_t\cdot\nabla u)\cdot u_t dx|\leq \|\nabla u\|_\i\|\sqrt\r
u_t\|_2^2,
\end{equation}
and
\begin{equation}
\begin{array}{ll}
\di \int\nabla(\l(\r)_t{\rm div}u)\cdot u_t dx&\di
=-\int\l(\r)_t{\rm div}u {\rm div} u_t dx\\
&\di\leq\f{\mu}{4}\int|{\rm div} u_t|^2 dx+C\|\l(\r)_t\|_4^2\|{\rm
div} u\|_4^2\di\leq\f{\mu}{4}\int|{\rm div} u_t|^2 dx+C.
\end{array}
\end{equation}
Substituting the above estimates into \eqref{u-tt1} and then
integrating with respect to $t$ over $[0,t]$ yield that
\begin{equation}\label{111}
\|\sqrt\r u_t\|_2^2(t)+\int_0^t\|\nabla u_t\|_2^2dt\leq \|\sqrt{\r^\d_0}
u^\d_t(0)\|_2^2+C\int_0^t(\|\nabla u\|_\i+1)\|\sqrt\r u_t\|^2_2 dt+C.
\end{equation}
By the compatibility condition \eqref{new-cc}, it holds that
$$
\sqrt{\r^\d_0} u^\d_t(0)=\f{\sqrt\r_0}{\sqrt{\r^\d_0}} g-\sqrt{\r^\d_0} u^\d_0\cdot\nabla u^\d_0,
$$
thus we have
$$
\|\sqrt{\r^\d_0} u^\d_t(0)\|_2^2\leq
\|\f{\sqrt\r_0}{\sqrt{\r^\d_0}} g\|_2^2+\|\r^\d_0\|_\i\|u^\d_0\|^2_\i\|\nabla u^\d_0\|_2^2\leq C,
$$
which, together with \eqref{111} and the Gronwall's inequality, yields that
\begin{equation}
\sup_{t\in[0,T]}\|\sqrt\r u_t\|_2^2(t)+\int_0^T\|\nabla
u_t\|_2^2dt\leq C.
\end{equation}
By \eqref{ut-HL}, for any $1\leq p<+\i,$
$$
\begin{array}{ll}
\di \int_0^T\|u_t\|_p^2dt&\di \leq
\int_0^T(\|(H,L)\|_p^2+\|u\|_\i^2\|\nabla u\|_p^2)dt\\
&\di \leq\int_0^T(\|\nabla(H,L)\|_2^2+\|u\|_\i^2\|\nabla
u\|_p^2)dt\leq C.
\end{array}
$$
Therefore, one can arrive at
$$
\int_0^T\|u_t\|_{H^1}^2dt\leq C.
$$
Thus the proof of Lemma \ref{nabla-ut-1} is completed. $\hfill\Box$

\begin{Lemma}\label{lemma-rho-tt}
It holds that
  \begin{equation}
\sup_{t\in[0,T]}\|(\r_t,P(\r)_t,\l(\r)_t)\|_{H^1}(t)+\int_0^T\|(\r_{tt},P(\r)_{tt},\l(\r)_{tt})\|_2^2
dt\leq C.
\end{equation}
\end{Lemma}
{\bf Proof:} From the continuity equation, it holds that $\r_t=-u\cdot \nabla
\r-\r{\rm div}u$ and
$\r_{tt}=-u_t\cdot\nabla\r-u\cdot\nabla\r_t-\r_t{\rm div}u-\r{\rm
div}u_t,$ and thus
\begin{equation}
\sup_{t\in[0,T]}\|\nabla\r_t\|_2(t)\leq
\sup_{t\in[0,T]}\Big[\|\nabla\r\|_4\|\nabla
u\|_4+\|u\|_\i\|\nabla^2\r\|_2+\|\r\|_\i\|\nabla^2u\|_2\Big]\leq C.
\end{equation}
and
\begin{equation}\label{rtt}
\begin{array}{ll}
\di \int_0^T\|\r_{tt}\|_2^2 dt&\di \leq
\int_0^T\Big[\|u_t\|_4^2\|\nabla\r\|_4^2+\|u\|_\i^2\|\nabla\r_t\|_2^2+\|\r_t\|_4^2\|\nabla
u\|_4^2+\|\r\|_\i^2\|\nabla u_t\|_2^2\Big]dt\\
&\di\leq C\int_0^T(\|u_t\|_{H^1}^2+1)dt\leq C.
\end{array}
\end{equation}
Similarly, we have
\begin{equation}
\sup_{t\in[0,T]}\|\nabla
(P(\r)_t,\l(\r)_t)\|_2(t)+\int_0^T\|(P(\r)_{tt},\l(\r)_{tt})\|_2^2
dt\leq C.
\end{equation}
Thus the proof of Lemma \ref{lemma-rho-tt} is completed. $\hfill\Box$

\begin{Lemma}\label{nabla-ut}
  It holds that
  \begin{equation}
  \begin{array}{ll}
    \di\sup_{t\in[0,T]}\Big[t\|u_t\|^2_{H^1}+t\|u\|^2_{H^3}+t\|(\r_{tt},P(\r)_{tt},\l(\r)_{tt})\|^2_2+\|(\r,P(\r))\|_{W^{2,q}}\Big]\\
    \di \qquad\qquad\qquad\qquad\qquad+\int_0^Tt\Big[\|\sqrt\r u_{tt}\|_2^2(t)+\|u_t\|_{H^2}^2(t)+\|u\|_{H^4}^2\Big] dt\leq C.
   \end{array} \end{equation}
\end{Lemma}
{\bf Proof:} Now multiplying the equation \eqref{u-tt} by $u_{tt}$ and then
integrating with respect to $x$ over $\mathbb{T}^2$ yield that
\begin{equation}\label{u-tt2}
\begin{array}{ll}
\di \|\sqrt\r u_{tt}\|_2^2(t)+\f12\f{d}{dt}\int\Big(\mu|\nabla
u_t|^2+(\mu+\l(\r))|{\rm div} u_t|^2\Big) dx=\f12\int\l(\r)_t|{\rm
div} u_t|^2 dx\\
\di -\int\big(\nabla P_t+\r_tu_t+\r_tu\cdot\nabla u+\r u\cdot \nabla
u_t+\r u_t\cdot \nabla u\big)\cdot u_{tt}dx+\int\nabla(\l(\r)_t{\rm
div} u) \cdot u_{tt}dx.
\end{array}
\end{equation}
Note that
$$
\begin{array}{ll}
\di \int\nabla(\l(\r)_t{\rm div} u) \cdot u_{tt}dx\di
=-\int\l(\r)_t{\rm
div} u{\rm div} u_{tt}dx\\
~~~\di =-\f{d}{dt}\int\l(\r)_t {\rm div} u{\rm div}
u_{t}dx+\int\big(\l(\r)_t|{\rm div} u|^2+\l(\r)_{tt}{\rm div} u{\rm
div} u_t\big)dx.
\end{array}
$$
Substituting the above identity into \eqref{u-tt2} yields that
\begin{equation}\label{u-tt3}
\begin{array}{ll}
\di \|\sqrt\r u_{tt}\|_2^2(t)+\f12\f{d}{dt}\int\Big(\mu|\nabla
u_t|^2+(\mu+\l(\r))|{\rm div} u_t|^2+\l(\r)_t {\rm div} u{\rm div}
u_{t}\Big) dx =\f32\int\l(\r)_t|{\rm div} u_t|^2 dx\\
\di~~~-\int\big(\nabla P_t+\r_tu_t+\r_tu\cdot\nabla u+\r u\cdot
\nabla u_t+\r u_t\cdot \nabla u\big)\cdot
u_{tt}dx+\int\l(\r)_{tt}{\rm div} u {\rm div} u_{t}dx.
\end{array}
\end{equation}
Note that $\l(\r)$ satisfies the transport equation
$\l(\r)_t=-u\cdot\nabla\l(\r)-\r\l^\prime(\r){\rm div}u,$ and then
it holds that
\begin{equation}\label{l1}
\begin{array}{ll}
\di |\f32\int\l(\r)_t|{\rm div} u_t|^2 dx|&\di =|-\f32\int u\cdot\nabla\l(\r)|{\rm div} u_t|^2 dx-\f32\int \r\l^\prime(\r){\rm div}u|{\rm div} u_t|^2dx|\\
&\di =|3\int\l(\r){\rm div} u_t u\cdot\nabla({\rm div} u_t)
dx+\f32\int(\l(\r)-\r\l^\prime(\r)){\rm div}u|{\rm div} u_t|^2 dx|\\
&\di \leq C\|\l(\r)u\|_\i\|{\rm div} u_t\|_2\|\nabla({\rm div}
u_t)\|_2+C\|\l(\r)-\r\l^\prime(\r)\|_\i\|\nabla u\|_\i\|{\rm div}
u_t\|_2^2\\
&\di \leq C\|{\rm div} u_t\|_2\|\nabla({\rm div} u_t)\|_2+C\|\nabla
u\|_\i\|{\rm div} u_t\|_2^2.
\end{array}
\end{equation}
It follows from \eqref{u-tt} that
$$
\mathcal{L}_\r u_t=\r u_{tt}+\r_tu_t+(\r u\cdot\nabla u)_t+\nabla
P(\r)_t+\nabla(\l(\r)_t{\rm div}u).
$$
Then the standard elliptic estimates show that
\begin{equation}\label{nabla-u2-t}
\begin{array}{ll}
\di \|\nabla^2u_t\|_{2}&\di\leq  C\Big[\|\sqrt\r\|_\i\|\sqrt\r
u_{tt}\|_2+\|\r_t\|_4\|u_t\|_4+\|\r_t\|_4\|u\|_\i\|\nabla
u\|_4+\|\r\|_\i\|u_t\|_4\|\nabla u\|_4\\
&\di~~~ +\|\r u\|_\i\|\nabla u_t\|_2+\|\nabla
P(\r)_t\|_2+\|\nabla\l(\r)_t\|_2\|{\rm
div}u\|_\i+\|\l(\r)_t\|_4\|\nabla^2 u\|_4\Big]\\
&\di \leq C\Big[\|\sqrt\r u_{tt}\|_2+\|u_t\|_4+1+\|\nabla
u_t\|_2+\|{\rm div}u\|_\i+\|\nabla^2 u\|_4\Big]\\
&\di \leq C\Big[\|\sqrt\r u_{tt}\|_2+\|u_t\|_4+1+\|\nabla
u_t\|_2+\|{\rm div}u\|_\i+\|\nabla^3 u\|_2\Big].
\end{array}
\end{equation}
Substituting \eqref{nabla-u2-t} into \eqref{l1} yields that
\begin{equation}\label{l2}
|\f32\int\l(\r)_t|{\rm div} u_t|^2 dx|\leq \f18\|\sqrt\r
u_{tt}\|_2^2+C(\|\nabla u\|_\i+1)\|\nabla
u_t\|_2^2+C\big(\|u_t\|_4^2+\|\nabla^3 u\|_2^2+\|\nabla u\|_\i^2).
\end{equation}
At the same time, it holds that
\begin{equation}
\begin{array}{ll}
\di -\int\nabla P(\r)_t\cdot u_{tt} dx &\di=\int P(\r)_t{\rm div} u_{tt}dx=\f{d}{dt}\int P(\r)_t{\rm div} u_{t}dx-\int P(\r)_{tt}{\rm div}
u_{t} dx\\
&\di \leq\f{d}{dt}\int P(\r)_t{\rm div}
u_{t}dx+\|P(\r)_{tt}\|_2^2+\|{\rm div} u_{t}\|_2^2,
\end{array}
\end{equation}
\begin{equation}
\begin{array}{ll}
\di -\int\r_tu_t\cdot u_{tt} dx&\di =\int \r_t\big(\f{|u_t|^2}{2}\big)_tdx=\f{d}{dt}\int \r_t\f{|u_t|^2}{2}dx-\int
\r_{tt}\f{|u_t|^2}{2}dx,
\end{array}
\end{equation}
while
\begin{equation}
\begin{array}{ll}
\di-\int
\r_{tt}\f{|u_t|^2}{2}dx&\di =\int {\rm div}(\r u)_t\f{|u_{t}|^2}{2}dx=\int(\r u)_t\cdot \nabla u_t\cdot u_t dx\\
&\di\leq \|\sqrt\r\|_\i\|\sqrt\r u_t\|_2\|u_t\|_4\|\nabla u_t\|_4+\|u\|_\i\|\r_t\|_4\|u_t\|_4\|\nabla u_t\|_2\\
&\di \leq C\Big[\|u_t\|_4\|\nabla u_t\|_4+\|u_t\|_4\|\nabla
u_t\|_2\Big] \leq C\Big[\|u_t\|_4\|\nabla^2 u_t\|_2+\|u_t\|_4\|\nabla
u_t\|_2\Big]\\
&\di \leq C\|u_t\|_4\Big[\|\sqrt\r u_{tt}\|_2+\|u_t\|_4+1+\|\nabla
u_t\|_2+\|{\rm div}u\|_\i+\|\nabla^3 u\|_2\Big]\\
&\di \leq \f18\|\sqrt\r u_{tt}\|_2^2+C\Big[\|u_t\|^2_4+1+\|\nabla
u_t\|^2_2+\|{\rm div}u\|^2_\i+\|\nabla^3 u\|^2_2\Big].
\end{array}
\end{equation}
Moreover, it follows that
\begin{equation}
\begin{array}{ll}
\di -\int\r_tu\cdot\nabla u\cdot u_{tt} dx\\
\di =-\f{d}{dt}\int \r_tu\cdot\nabla u\cdot u_{t}dx+\int
\r_{tt}u\cdot\nabla u\cdot u_{t}dx\int\r_tu_t\cdot\nabla u\cdot
u_{t}dx+\int
\r_tu\cdot\nabla u_t\cdot u_{t}dx\\
 \di =-\f{d}{dt}\int \r_tu\cdot\nabla u\cdot u_{t}dx+\|\r_{tt}\|_2\|u\|_\i\|\nabla u\|_4\|u_{t}\|_4\\
\di~~~~~~~~~~~~~~+\|\r_t\|_4\|u_t\|_4^2\|\nabla
u\|_4+\|\r_t\|_4\|u\|_\i\|\nabla u_t\|_2\|u_{t}\|_4\\
 \di
\leq -\f{d}{dt}\int \r_tu\cdot\nabla u\cdot
u_{t}dx+C\Big[\|\r_{tt}\|_2\|u_{t}\|_4+\|u_t\|_4^2+\|\nabla
u_t\|_2\|u_{t}\|_4\Big]\\
\di\leq -\f{d}{dt}\int \r_tu\cdot\nabla u\cdot
u_{t}dx+C\Big[\|\r_{tt}\|_2^2+\|u_t\|_4^2+\|\nabla u_t\|_2^2\Big],
\end{array}
\end{equation}
\begin{equation}
\begin{array}{ll}
\di -\int\r u\cdot\nabla u_t\cdot u_{tt} dx &\di\leq \|\sqrt\r u_{tt}\|_2\|\sqrt\r u\|_\i\|\nabla u_t\|_2 \leq\f18\|\sqrt\r u_{tt}\|_2^2+C\|\nabla u_{t}\|_2^2,
\end{array}
\end{equation}
\begin{equation}
\begin{array}{ll}
\di -\int\r u_t\cdot\nabla u\cdot u_{tt} dx&\di\|\sqrt\r u_{tt}\|_2\|\sqrt\r\|_\i\|\nabla u\|_4\|u_t\|_4 \leq\f18\|\sqrt\r u_{tt}\|_2^2+C\| u_{t}\|_4^2,
\end{array}
\end{equation}
and
\begin{equation}
\begin{array}{ll}
\di -\int\l(\r)_{tt}{\rm div}u {\rm div} u_{t} dx&\di\leq \|\l(\r)_{tt}\|_2\|\nabla u\|_\i\|{\rm div}u_t\|_2 \leq \f12\Big[\|\l(\r)_{tt}\|_2^2+\|\nabla u\|^2_\i\|{\rm
div}u_t\|_2^2\Big].
\end{array}
\end{equation}
Collecting all the above estimates and substituting them into
\eqref{u-tt3} yield that
\begin{equation}\label{du2}
\begin{array}{ll}
\di \f12\|\sqrt\r u_{tt}\|_2^2(t)+\f{d}{dt} G(t)\\
\di \leq
C\Big[\|(\r_{tt},P(\r)_{tt},\l(\r)_{tt})\|_2^2+\|u_t\|_4^2+\|\nabla^3
u\|_2^2+(\|\nabla u\|_\i^2+1)(\|\nabla u_t\|_2^2+1)\Big]
\end{array}
\end{equation}
where
\begin{equation}
G(t)=\int\Big(\mu|\nabla u_t|^2+(\mu+\l(\r))|{\rm div}
u_t|^2+\l(\r)_t {\rm div} u{\rm div} u_{t}-P(\r)_t{\rm div} u_t
+\r_t\f{|u_t|^2}{2}+\r_t u\cdot \nabla u\cdot u_t\Big) dx.
\end{equation}
Note that
$$
\begin{array}{ll}
\di |\int\l(\r)_t {\rm div} u{\rm div} u_{t} dx|&\di \leq
\|\l(\r)_t\|_4\| {\rm div} u\|_4\|{\rm div} u_{t}\|_2\\
&\di\leq \f{\mu}{8}\|\nabla u_t\|_2^2+C\|\l(\r)_t\|^2_4\| {\rm div}
u\|_4^2\leq \f{\mu}{8}\|\nabla u_t\|_2^2+C,
\end{array}
$$
$$
|-\int P(\r)_t{\rm div} u_tdx|\leq \f{\mu}{8}\|\nabla
u_t\|_2^2+C\|P(\r)_t\|_2^2\leq \f{\mu}{8}\|\nabla u_t\|_2^2+C,
$$
$$
\begin{array}{ll}
\di|\int\r_t\f{|u_t|^2}{2} dx|&\di =|\int{\rm div}(\r
u)\f{|u_t|^2}{2} dx|=|\int\r u\cdot \nabla u_t\cdot u_t dx|\\
&\di \leq \|\sqrt\r u_t\|_2\|\sqrt\r u\|_\i\|\nabla u_t\|_2\leq
\f{\mu}{8}\|\nabla u_t\|_2^2+C,
\end{array}
$$
and
$$
\begin{array}{ll}
\di|\int\r_t u\cdot \nabla u\cdot u_tdx|&\di =|\int{\rm div}(\r
u) (u\cdot \nabla u\cdot u_t) dx|=|\int\r u\cdot \nabla (u\cdot \nabla u\cdot u_t)dx|\\
&\di \leq \|\sqrt\r u_t\|_2\|\sqrt\r u\|_\i\big(\|\nabla
u\|_4^2+\|u\|_\i\|\nabla^2 u\|_2\big)+\|\r |u|^2\|_\i\|\nabla
u_t\|_2\|\nabla u\|_2\\
&\di \leq \f{\mu}{8}\|\nabla u_t\|_2^2+C.
\end{array}
$$
Therefore, it holds that
\begin{equation}\label{G}
C_1(\|\nabla u_t\|_2^2-1)\leq G(t)\leq C(\|\nabla u_t\|_2^2+1),
\end{equation}
for some positive constants $C,C_1$.

Now from \eqref{du2}, we can arrive at
\begin{equation}\label{du1}
\begin{array}{ll}
\di \f12\|\sqrt\r u_{tt}\|_2^2(t)+\f{d}{dt} G(t)\\
\di \leq
C\Big[\|(\r_{tt},P(\r)_{tt},\l(\r)_{tt})\|_2^2+\|u_t\|_4^2+\|\nabla^3
u\|_2^2+(\|\nabla u\|_\i^2+1)(G(t)+1)\Big].
\end{array}
\end{equation}
Multiplying the above inequality by $t$ and then integrating the resulting inequality with respect to $t$ over the interval $[\tau, t_1]$ with $\tau, t_1\in[0,T]$ give that
\begin{equation}\label{tau-k}
\begin{array}{ll}
\di \int_{\tau}^{t_1}t\|\sqrt\r u_{tt}\|_2^2(t)dt+t_1G(t_1)\leq C\tau G(\tau)+C\int_{\tau}^{t_1}\Big[(\|\nabla u\|_\i^2+1)(tG(t)+1)\Big]dt\\
\di \qquad\qquad\qquad\qquad\qquad+
C\int_{\tau}^{t_1}\Big[\|(\r_{tt},P(\r)_{tt},\l(\r)_{tt})\|_2^2+\|u_t\|_4^2+\|\nabla^3
u\|_2^2+G(t)\Big]dt.
\end{array}
\end{equation}
It follows from Lemma \ref{nabla-ut-1} and \eqref{G} that $G(t)\in L^1(0,T)$. Thus, due to \cite{CCK}, there exists a subsequence $\tau_k$ such that
\begin{equation}
\tau_k\rightarrow 0,\qquad \tau_kG(\tau_k)\rightarrow 0,\qquad {\rm as}~~ k\rightarrow +\i.
\end{equation}
Taking $\tau=\tau_k$ in \eqref{tau-k}, then $k\rightarrow+\i$ and using the Gronwall's  inequality, one gets that
\begin{equation}
\sup_{t\in[0,T]}\big[t\|\nabla u_t\|_2^2(t)\big]+\int_0^Tt\|\sqrt\r
u_{tt}\|_2^2(t)dt\leq C.
\end{equation}
Note that \eqref{nabla-u2-t} implies that
\begin{equation}
 \sup_{t\in[0,T]}\big[t\|(\r_{tt},P(\r)_{tt},\l(\r)_{tt})\|_2^2(t)\big]+ \int _0^Tt\|\nabla^2
u_{t}\|_2^2(t)dt\leq C.
\end{equation}
It follows from \eqref{ut-HL} that
$$
\nabla(L,H)^t=\nabla u_t-\nabla (u\cdot\nabla u).
$$
Consequently, it holds that
\begin{equation}
\sup_{t\in[0,T]}\big[t\|\nabla(L,H)^t\|^2_2(t)\big]\leq C,
\end{equation}
which, together with \eqref{fact5}, implies that
\begin{equation}
\sup_{t\in[0,T]}\big[t\|(L,H)^t\|^2_2(t)\big]\leq C\sup_{t\in[0,T]}\big[t\|\nabla(L,H)^t\|^2_2(t)\big]\leq C.
\end{equation}
Therefore, it follows that
\begin{equation}
\sup_{t\in[0,T]}\big[t\|u_t\|^2_2(t)\big]\leq C\sup_{t\in[0,T]}\big[t\|(L,H)^t\|^2_2(t)+t\|u\cdot\nabla u\|_2^2(t)\big]\leq C.
\end{equation}
So one can infer further that
\begin{equation}
\sup_{t\in[0,T]}\big[t\|u_t\|^2_{H^1}(t)\big]+\int_0^Tt\|u_t\|_{H^2}^2(t)dt\leq C.
\end{equation}
Applying $\partial_{x_jx_k}$, $j,k=1,2$, to $\eqref{CNS}_1$ gives
$$
\begin{array}{ll}
\di (\r_{x_jx_k})_t+u\cdot \nabla(\r_{x_jx_k})+u_{x_jx_k}\cdot\nabla\r+u_{x_j}\cdot\nabla\r_{x_k}+u_{x_k}\cdot\nabla\r_{x_j}\\
\di\qquad~~~+\r_{x_jx_k}{\rm div} u+\r_{x_j}({\rm div} u)_{x_k}+\r_{x_k}({\rm div} u)_{x_j}+\r({\rm div} u)_{x_jx_k}=0.
\end{array}
$$
Multiplying the above equation by $q|\nabla^2\r|^{q-2}\r_{x_jx_k}$ with $q>2$ given in Theorem \ref{theorem2} and summing over $j,k=1,2$ give that
$$
\begin{array}{ll}
\di (|\nabla^2\r|^q)_t+{\rm div}(u|\nabla^2\r|^q)+(q-1)|\nabla^2\r|^q{\rm div}u
+q|\nabla^2\r|^{q-2}\r_{x_jx_k}\Big[u_{x_jx_k}\cdot\nabla\r\\
\di +u_{x_j}\cdot\nabla\r_{x_k}+u_{x_k}\cdot\nabla\r_{x_j}+\r_{x_j}({\rm div} u)_{x_k}+\r_{x_k}({\rm div} u)_{x_j}+\r({\rm div} u)_{x_jx_k}\Big]=0.
\end{array}
$$
Integrating the above equality with respect to $x$ over $\mathbb{T}^2$ leads to that
$$
\f{d}{dt}\|\nabla^2\r\|_q^q\leq (q-1)\|{\rm div}
u\|_\i\|\nabla^2\r\|_q^q+Cq\|\nabla^2\r\|_q^{q-1}\Big[\|\nabla\r\|_{2q}\|\nabla^2
u\|_{2q}+\|\nabla
u\|_\i\|\nabla^2\r\|_q+\|\r\|_\i\|\nabla^3u\|_q\Big].
$$
Thus one can get
\begin{equation}\label{r1}
\begin{array}{ll}
\di\f{d}{dt}\|\nabla^2\r\|_q&\di \leq C\Big[\|\nabla u\|_\i\|\nabla^2\r\|_q+\|\nabla\r\|_{2q}\|\nabla^2 u\|_{2q}+\|\r\|_\i\|\nabla^3u\|_q\Big]\\
&\di \leq C\Big[\|\nabla u\|_\i\|\nabla^2\r\|_q+\|\nabla^2
u\|_{W^{1,q}}\Big],
\end{array}
\end{equation}
where $q>2$.  Similarly, one can obtain
\begin{equation}\label{P1}
\di\f{d}{dt}\|\nabla^2P\|_q \leq C\Big[\|\nabla u\|_\i\|\nabla^2P\|_q+\|\nabla^2 u\|_{W^{1,q}}\Big].
\end{equation}
Apply $\partial_{x_i}$ with $i=1,2$ to the elliptic system
 $
\mathcal{L}_\r u=\r u_t+\r u\cdot\nabla u+\nabla P(\r)
$
to get
$$
\mathcal{L}_\r u_{x_i}=-\nabla(\l(\r)_{x_i}{\rm div} u)+\r_{x_i} u_t+\r u_{x_i t}+\r_{x_i} u\cdot\nabla u+\r u_{x_i}\cdot\nabla u+\r u\cdot\nabla u_{x_i}+\nabla P(\r)_{x_i}:=\Psi.
$$
Then the standard elliptic regularity estimates imply that
\begin{equation}\label{w2q}
\begin{array}{ll}
\di \|\nabla u\|_{W^{2,q}}&\di\leq C\Big[\|\nabla u\|_q+\|\Psi\|_q\Big]\\
&\di \leq C\Big[1+(\|\nabla u\|_\i+1)\|(\nabla^2\r,\nabla^2P)\|_q+\|\nabla u\|_{W^{1,q}}+\|u_t\|_{W^{1,q}}\Big]\\
&\di \leq C\Big[1+(\|\nabla u\|_\i+1)\|(\nabla^2\r,\nabla^2P)\|_q+\|u\|_{H^3}+\|u_t\|_{H^1}+\|\nabla u_t\|_{q}\Big].
\end{array}
\end{equation}
Thus it follows from \eqref{r1}, \eqref{P1} and \eqref{w2q} that
\begin{equation}\label{r2}
\f{d}{dt}\|(\nabla^2\r,\nabla^2P)\|_q\leq C\Big[1+(\|\nabla u\|_\i+1)\|(\nabla^2\r,\nabla^2P)\|_q+\|u\|_{H^3}+\|u_t\|_{H^1}+\|\nabla u_t\|_{q}\Big].
\end{equation}
Note that Lemma \ref{lemma2} implies that
$$
\begin{array}{ll}
\di\int_0^T\|\nabla u_t\|_{q}(t)dt&\di \leq C \int_0^T\|\nabla^2
u_t\|_2(t)dt \leq C \sup_{t\in[0,T]}\big[\sqrt t\|\nabla^2
u_t\|_2(t)\big]\int_0^T t^{-\f12}dt\leq C.
\end{array}
$$
Therefore, it follows from \eqref{r2} and the Gronwall's inequality that
\begin{equation}\label{2dp}
\begin{array}{ll}
\di \|(\nabla^2\r,\nabla^2P(\r))\|_{q}(t)\leq \Big(\|(\nabla^2\r_0,\nabla^2P(\r_0))\|_{q}\\
\di \qquad\qquad+C\int_0^t(1+\|u\|_{H^3}+\|u_t\|_{H^1}+\|\nabla
u_t\|_{q})ds\Big)e^{\di C\int_0^t(\|\nabla u\|_\i(s)+1)ds} \leq C,
\end{array}
\end{equation}
which then gives
\begin{equation}
\sup_{t\in[0,T]}\|(\r,P(\r))\|_{W^{2,q}(\mathbb{T}^2)}\leq C.
\end{equation}
So the proof of Lemma \ref{nabla-ut} is completed. $\hfill\Box$
\begin{Lemma}\label{lemma4.6}
  It holds that for any $0<\tau\leq T$,
  \begin{equation}
    \sup_{t\in[0,T]}\big[t^2\|\sqrt\r u_{tt}\|_2^2(t)+t^2\|u_t\|_{H^2}^2+t^2\|u\|_{W^{3,q}}^2\big]+\int_0^Tt^2\|\nabla u_{tt}\|_2^2(t) dt\leq C.
  \end{equation}
\end{Lemma}
{\bf Proof:} Applying $\partial_t$ to the equation \eqref{u-tt} gives that
\begin{equation}
\label{u-ttt}
\begin{array}{ll}
&\di \r u_{ttt}+\r u\cdot \nabla u_{tt}-\mathcal{L}_\r u_{tt}=-\nabla p_{tt}-\r_{tt}(u_t+u\cdot\nabla u)-2\r_t(u_{tt}+u_t\cdot\nabla u+u\cdot\nabla u_t)\\
&\di-2\r u_t\cdot\nabla u_t-\r u_{tt}\cdot \nabla u+2\nabla((\l(\r))_t{\rm div} u_t)+\nabla(\l(\r)_{tt}{\rm div}u).
\end{array}
\end{equation}
Multiplying the equation \eqref{u-ttt} by $u_{tt}$ and integrating the resulting equation with respect to $x$ over $\mathbb{T}^2$ yield that
$$
\begin{array}{ll}
\di\f12\f{d}{dt}\int \r|u_{tt}|^2dx+\int\mu|\nabla u_{tt}|^2+(\mu+\l(\r))({\rm div} u_{tt})^2 dx=\int p_{tt}{\rm div}u_{tt} dx\\
\di\qquad -\int\r_{tt}(u_t+u\cdot\nabla u)\cdot u_{tt} dx-2\int\r_t(u_{tt}+u_t\cdot\nabla u+u\cdot\nabla u_t)\cdot u_{tt} dx-2\int\r u_t\cdot\nabla u_t \cdot u_{tt}dx\\
\di\qquad -\int \r u_{tt}\cdot \nabla u\cdot u_{tt}dx -2\int \l(\r)_t{\rm div}u_t {\rm div}u_{tt}dx-\int\l(\r)_{tt}{\rm div}u{\rm div} u_{tt} dx.
\end{array}
$$
Multiply the above equality by $t^2$ to get that
\begin{equation}\label{he}
\begin{array}{ll}
\di\f12\f{d}{dt}\Big(t^2\int \r|u_{tt}|^2dx\Big)-t\int \r|u_{tt}|^2dx +t^2\int\mu|\nabla u_{tt}|^2+(\mu+\l(\r))({\rm div} u_{tt})^2 dx=t^2\int P_{tt}{\rm div}u_{tt} dx\\
\di-t^2\int\r_{tt}(u_t+u\cdot\nabla u)\cdot u_{tt} dx-2t^2\int\r_t(u_{tt}+u_t\cdot\nabla u+u\cdot\nabla u_t)\cdot u_{tt} dx-2t^2\int\r u_t\cdot\nabla u_t \cdot u_{tt}dx\\
\di-t^2\int \r u_{tt}\cdot \nabla u\cdot u_{tt}dx -2t^2\int \l(\r)_t{\rm div}u_t {\rm div}u_{tt}dx-t^2\int\l(\r)_{tt}{\rm div}u{\rm div} u_{tt} dx:=\sum_{i=1}^7I_i.
\end{array}
\end{equation}
Clearly,
$$
|I_1|\leq \a t^2\|{\rm div} u_{tt}\|_2^2+C_\a t^2\|P_{tt}\|_2^2.
$$
Now we estimate $I_2$, which is a little more delicate due to the absence of estimates for $u_{tt}$.  First, rewrite $I_2$ as
$$
\begin{array}{ll}
\di I_2&\di=t^2\int{\rm div}(\r u)_{t}(L,H)^t\cdot u_{tt} dx =-t^2\int(\r u)_{t}\cdot\nabla\big[(L,H)^t\cdot u_{tt}\big] dx\\
&\di =-t^2\int\r u_t \cdot\nabla\big[(L,H)^t\cdot u_{tt}\big]dx-t^2\int\r_tu\cdot\nabla u_{tt}\cdot (L,H)^t dx-t^2\int\r_tu\cdot\nabla(L,H)^t\cdot u_{tt} dx\\
&\di =-t^2\int\r u_t \cdot\nabla\big[(L,H)^t\cdot u_{tt}\big]dx-t^2\int\r_tu\cdot\nabla u_{tt}\cdot (L,H)^t dx\\
&\di\qquad\qquad\qquad\qquad\qquad -t^2\int\r
u\cdot\nabla\big[u\cdot\nabla(L,H)^t\cdot u_{tt}\big]
dx:=I_{21}+I_{22}+I_{23}
\end{array}
$$
where the superscript $^t$ means the transpose of the vector
$(L,H)$.

 Now, direct estimates yields that
\begin{equation}
\begin{array}{ll}
 |I_{21}|&\di\leq t^2\|\sqrt\r u_{tt}\|_2\|\sqrt\r\|_\i\|u_t\|_\i\|\nabla(L,H)^t\|_2+t^2\|\r\|_\i\|\nabla u_{tt}\|_2\|u_t\|_4\|(L,H)^t\|_4\\
 &\di \leq Ct^2\Big[\|\sqrt\r u_{tt}\|_2\|u_t\|_{H^2}\|\nabla(L,H)^t\|_2+\|\nabla u_{tt}\|_2\|u_t\|_{H^1}\|\nabla(L,H)^t\|_2\Big]\\
  &\di \leq \a t^2\|\nabla u_{tt}\|_2^2+C_\a t^2\Big[\|\sqrt\r u_{tt}\|_2^2\|\nabla(L,H)^t\|^2_2+\|u_t\|_{H^2}^2+\|u_t\|_{H^1}^2\|\nabla(L,H)^t\|^2_2\Big]\\
  &\di \leq \a t^2\|\nabla u_{tt}\|_2^2+C_\a \Big[t^2\|\sqrt\r u_{tt}\|_2^2\|\nabla(L,H)^t\|^2_2+t^2\|u_t\|_{H^2}^2+\|\nabla(L,H)^t\|^2_2\Big]
\end{array}
 \end{equation}
where in the last inequality one has used Lemma \ref{nabla-ut}.

Similarly, one can obtain
\begin{equation}
\begin{array}{ll}
  |I_{22}|&\di \leq t^2\|\nabla u_{tt}\|_2\|u\|_\i\|\r_t\|_4\|(L,H)^t\|_4\\
  &\di \leq \a t^2\|\nabla u_{tt}\|_2^2+C_\a t^2\|\r_t\|_{H^1}^2\|\nabla(L,H)^t\|_2^2\leq \a t^2\|\nabla u_{tt}\|_2^2+C_\a \|\nabla(L,H)^t\|_2^2
\end{array}
\end{equation}
and
\begin{equation}
\begin{array}{ll}
  |I_{23}|&\di \leq t^2\Big[\|\sqrt\r u_{tt}\|_2\|\sqrt\r u\|_\i\|\nabla u\|_\i\|\nabla(L,H)^t\|_2+\|\nabla u_{tt}\|_2\|\r u^2\|_\i\|\nabla(L,H)^t\|_2\\
  &\di\qquad\qquad\qquad+\|\sqrt\r u_{tt}\|_2\|\sqrt\r u^2\|_\i(\|\nabla^2 u_t\|_2+\|u\|_\i\|\nabla^3u\|_2+\|\nabla u\|_\i\|\nabla^2u\|_2)\Big]\\
  &\di \leq \a t^2\|\nabla u_{tt}\|_2^2+C_\a \|\nabla(L,H)^t\|_2^2\\
  &\di \qquad+C\Big[t^2\|\sqrt\r u_{tt}\|_2^2(t\|\nabla u\|_\i^2+1)+t^2\|\nabla^2 u_t\|_2^2+t^2\|u\|_{H^3}^2+\|\nabla(L,H)^t\|_2^2\Big].
\end{array}
\end{equation}
Continuing, using the lemmas obtained so far, one can get that
\begin{equation}
\begin{array}{ll}
|I_3|&\di \leq t^2\Big[\|\sqrt\r u\|_\i\|\nabla u_{tt}\|_2\|\sqrt\r u_{tt}\|_2+\|\r u\|_\i\|\nabla u_{tt}\|_2(\|u_t\cdot\nabla u\|_2+\|u\cdot \nabla u_t\|_2)\\
&\di \qquad\quad+\|\sqrt\r u\|_\i\|\sqrt\r u_{tt}\|_2(\|\nabla(u_t\cdot\nabla u)\|_2+\|\nabla(u\cdot\nabla u_t)\|_2)\Big]\\
&\di \leq t^2\Big[\|\sqrt\r u\|_\i\|\nabla u_{tt}\|_2\|\sqrt\r u_{tt}\|_2+\|\r u\|_\i\|\nabla u_{tt}\|_2(\|u_t\|_4\|\nabla u\|_4+\|u\|_\i\| \nabla u_t\|_2)\\
&\di \qquad\quad+\|\sqrt\r u\|_\i\|\sqrt\r u_{tt}\|_2(\|u_t\|_4\|\nabla^2 u\|_4+\|\nabla  u_t\|_4\|\nabla u\|_4+\|u\|_\i\|\nabla^2u_t\|_2)\Big]\\
&\di \leq \a t^2\|\nabla u_{tt}\|_2^2+C_\a t^2\Big[\|\sqrt\r u_{tt}\|_2^2+\|u_t\|_{H^1}^2(\|\nabla^3 u\|^2_2+1)+\|u_t\|_{H^2}^2\Big],
\end{array}
\end{equation}
\begin{equation}
\begin{array}{ll}
\di |I_4|&\di \leq t^2\|\sqrt\r u_{tt}\|_2\|\sqrt\r\|_\i\|u_t\|_4\|\nabla u_t\|_4\\
&\di\leq Ct^2\|\sqrt\r u_{tt}\|_2\|u_t\|_{H^1}\|\nabla^2 u_t\|_2
\leq Ct^2\|\sqrt\r u_{tt}\|_2^2+Ct^2\|u_t\|^2_{H^1}\|\nabla^2
u_t\|^2_2,
\end{array}
\end{equation}
\begin{equation}
|I_5|\leq Ct^2\|\sqrt\r u_{tt}\|^2_2\|\nabla u\|_\i,
\end{equation}
and
\begin{equation}
\begin{array}{ll}
\di|I_6|&\di \leq t^2\|{\rm div}u_{tt}\|_2\|\l(\r)_t\|_4\|\nabla
u_t\|_4\leq \a t^2\|{\rm div} u_{tt}\|_2^2+C_\a t^2\|\nabla^2
u_t\|_2^2,
\end{array}
\end{equation}
\begin{equation}
\begin{array}{ll}
\di |I_7|&\di \leq t^2\|{\rm div}u_{tt}\|_2\|\l(\r)_{tt}\|_2\|\nabla
u\|_\i\leq \a t^2\|{\rm div} u_{tt}\|_2^2+C_\a
t^2\|\l(\r)_{tt}\|_2^2\|u\|_{H^3}^2.
\end{array}
\end{equation}
Substituting the above estimates on $I_i~(i=1,2,\cdots,7)$ into \eqref{he} and then integrating the resulting inequality with respect $t$ over $[\tau,t_1]$ with $\tau,t_1\in[0,T]$ give that
\begin{equation}\label{rho-u-tt}
t_1^2\|\sqrt\r u_{tt}(t_1)\|_2^2+\int_{\tau}^{t_1}t^2\|\nabla u_{tt}\|_2^2 dt\leq C+C\tau^2\|\sqrt\r u_{tt}(\tau)\|_2^2.
\end{equation}
Since $t\sqrt\r u_{tt}\in L^2([0,T]\times\mathbb{T}^2)$ due to Lemma \ref{nabla-ut}, there exists a subsequence $\tau_k$ such that
\begin{equation}
\tau_k\rightarrow 0,\qquad \tau_k^2\|\sqrt\r u_{tt}(\tau_k)\|_2^2\rightarrow 0,\qquad {\rm as}~~ k\rightarrow +\i.
\end{equation}
Letting $\tau=\tau_k$ in \eqref{rho-u-tt} and $k\rightarrow +\i$, one gets that
\begin{equation}\label{rho-u-tt-1}
 t^2\|\sqrt\r u_{tt}(t)\|_2^2+\int_{0}^{t}s^2\|\nabla u_{tt}(s)\|_2^2 dt\leq C.
\end{equation}
By, \eqref{nabla-u2-t}, it holds that
\begin{equation}
\sup_{t\in[0,T]}\big[t^2\|\nabla^2 u_t\|_2^2(t)\big]\leq C\sup_{t\in[0,T]}\Big[t^2\|\sqrt\r u_{tt}\|_2^2(t)+t^2\|u_t\|_{H^1}^2+t^2\|u\|_{H^3}^2+1\Big]\leq C.
\end{equation}
Finally, by \eqref{w2q}, we can obtain
\begin{equation}
\sup_{t\in[0,T]}\big[t^2\|\nabla u\|_{W^{2,q}}^2(t)\big]\leq C\sup_{t\in[0,T]}\Big[t^2\|u\|_{H^3}^2+t^2\|u_t\|_{H^2}^2+1\Big]\leq C.
\end{equation}
So the proof of Lemma \ref{lemma4.6} is completed. $\hfill\Box$

\section{The proof of Theorem \ref{theorem2}}
With the uniform-in-$\d$ bounds of the solution $(\r^\d,u^\d)$ in
Lemmas \ref{lemma-ee}-\ref{lemma-F-o} and Lemma
\ref{lemma4.1}-\ref{lemma4.6}, one can prove the convergence of the
sequence $(\r^\d,u^\d)$ to a limit $(\r,u)$ satisfying the same
bounds as $(\r^\d,u^\d)$ as $\d$ tends to zero and the limit
$(\r,u)$ is a unique solution to the original problem
\eqref{CNS}-\eqref{initial-v}. The details are omitted for brevity and
one can refer to Cho-Kim \cite{CK} for the routine proofs. In the following, we will
show that $(\r,u)$ satisfy the bounds in Theorem \ref{theorem2} and
$(\r,u)$ is a classical solution to \eqref{CNS}. Since
$u\in L^2(0,T;H^3(\mathbb{T}^2))$ and $u_t\in
L^2(0,T;H^1(\mathbb{T}^2))$, so the Sobolev's
embedding theorem implies that
$$
u\in C([0,T];H^2(\mathbb{T}^2))\hookrightarrow C([0,T]\times\mathbb{T}^2).
$$
Then it follows from $(\r,P(\r))\in L^\i(0,T;W^{2,q}(\mathbb{T}^2))$ and $(\r,P(\r))_t\in L^\i(0,T;H^1(\mathbb{T}^2))$ that $(\r,P(\r))\in C([0,T];W^{1,q}(\mathbb{T}^2))\cap C([0,T];W^{2,q}(\mathbb{T}^2)-weak)$. This and  \eqref{2dp} then
imply that
$$
(\r,P(\r))\in C([0,T];W^{2,q}(\mathbb{T}^2)).
$$
Since for any $\tau\in(0,T)$,
$$
(\nabla u,\nabla^2 u)\in L^\i(\tau,T;W^{1,q}(\mathbb{T}^2)),\qquad (\nabla u_t,\nabla^2 u_t)\in  L^\i(\tau,T;L^2(\mathbb{T}^2)).
$$
Therefore,
$$
(\nabla u,\nabla^2u)\in C([\tau,T]\times\mathbb{T}^2),
$$
Due to the fact that
$$
\nabla (\r,P(\r))\in C([0,T];W^{1,q}(\mathbb{T}^2))\hookrightarrow C([0,T]\times\mathbb{T}^2)
$$
and the continuity equation $\eqref{CNS}_1$, it holds that
$$
\r_t=u\cdot\nabla \r+\r{\rm div}u\in C([\tau,T]\times\mathbb{T}^2).
$$
It follows from the momentum equation $\eqref{CNS}_2$ that
$$
\begin{array}{ll}
  (\r u)_t&\di =\mathcal{L}_{\r}u-{\rm div}(\r u\otimes u)-\nabla P(\r)\\
  &\di =\mu \Delta u+(\mu+\l(\r))\nabla({\rm div} u)+({\rm div}u)\nabla\l(\r)+\r u\cdot\nabla u+\r u{\rm div} u+(u\cdot\nabla\r) u-\nabla P(\r)\\
  &\di \in C([\tau,T]\times\mathbb{T}^2).
\end{array}
$$
Thus we completed the proof of Theorem \ref{theorem2}.

\section{The proof of Theorem \ref{theorem}}
Based on Theorem  \ref{theorem2}, one can prove Theorem \ref{theorem} easily as follows. Since
$$
\r_0\in H^3(\mathbb{T}^2)\hookrightarrow W^{2,q}(\mathbb{T}^2)
$$
for any $2<q<+\i$, it follows that under the conditions of Theorem \ref{theorem}, Theorem \ref{theorem2} holds for any $2<q<+\i$.
Thus, we need only to prove the higher order regularity presented  in Theorem \ref{theorem}.

\begin{Lemma}\label{new}
  It holds that
$$
  \begin{array}{ll}
    \di\sup_{t\in[0,T]}\Big[\|\sqrt\r \nabla^3u\|_{2}(t)+\|(\r,P(\r),\l(\r))\|_{H^3}(t)\Big]+\int_0^T\|u\|_{H^4}^2 dt\leq C.
   \end{array} $$
\end{Lemma}
\noindent{\bf Proof:}
Applying $\partial_{x_jx_k}$, $j,k=1,2$,  to $\eqref{moment-1}$ yields that
\begin{equation}\label{Jan-19-2}
 \begin{array}{ll}
\di \r u_{x_jx_kt}+\r u\cdot\nabla u_{x_jx_k}+\r_{x_jx_k}u_t+\r_{x_j}u_{x_k t}+\r_{x_k}u_{x_jt}+\r_{x_jx_k}u\cdot\nabla u+\r u_{x_jx_k}\cdot \nabla u\\
\di +\r_{x_j}u_{x_k}\cdot\nabla u+\r_{x_j}u\cdot\nabla u_{x_k}+\r_{x_k}u_{x_j}\cdot\nabla u+\r_{x_k}u\cdot\nabla u_{x_j}
+\r u_{x_j}\cdot\nabla u_{x_k}+\r u_{x_k}\cdot\nabla u_{x_j}\\
\di+\nabla P(\r)_{x_jx_k}=\mu\Delta u_{x_jx_k}+\nabla((\mu+\l(\r)){\rm div} u)_{x_jx_k}.
\end{array}
\end{equation}
Then multiplying $\eqref{Jan-19-2}$ by $\Delta u_{x_jx_k}$ and integrating with respect to $x$ over $\mathbb{T}^2$ imply that
\begin{equation}\label{u3}
\begin{array}{ll}
\di\int \big[\mu|\Delta u_{x_jx_k}|^2+\nabla((\mu+\l(\r)){\rm div} u)_{x_jx_k}\cdot\Delta u_{x_jx_k}\big]dx=\int\big(\r u_{x_jx_kt}+\r u\cdot\nabla u_{x_jx_k}\big)\cdot \Delta u_{x_jx_k} dx\\
\di +\int\big[ \r_{x_jx_k}u_t+\r_{x_j}u_{x_k t}+\r_{x_k}u_{x_jt}+\r_{x_jx_k}u\cdot\nabla u+\r u_{x_jx_k}\cdot \nabla u+\r_{x_j}u_{x_k}\cdot\nabla u+\r_{x_j}u\cdot\nabla u_{x_k}\\
\di\qquad\quad +\r_{x_k}u_{x_j}\cdot\nabla u+\r_{x_k}u\cdot\nabla u_{x_j}
+\r u_{x_j}\cdot\nabla u_{x_k}+\r u_{x_k}\cdot\nabla u_{x_j}+\nabla P(\r)_{x_jx_k}\big]\cdot \Delta u_{x_jx_k}dx.
\end{array}
\end{equation}
Integrations by part several times yield
\begin{equation}\label{u31}
\begin{array}{ll}
\di \int\big(\r u_{x_jx_kt}+\r u\cdot\nabla u_{x_jx_k}\big)\cdot \Delta u_{x_jx_k} dx\\
\di =-\int\Big[\r \Big(\f{|\nabla u_{x_jx_k}|^2}{2}\Big)_t+\r u\cdot \nabla\Big(\f{|\nabla u_{x_jx_k}|^2}{2}\Big)+\sum_{i=1}^2\r_{x_i}u_{x_jx_kt}\cdot u_{x_ix_jx_k}\\
\di\qquad\qquad\qquad+\sum_{i=1}^2\r_{x_i}u\cdot\nabla u_{x_jx_k}\cdot u_{x_ix_jx_k}+\sum_{i=1}^2\r u_{x_i}\cdot\nabla u_{x_jx_k}\cdot u_{x_ix_jx_k}\Big]dx\\
\di =-\f{d}{dt}\int\r \f{|\nabla u_{x_jx_k}|^2}{2}dx +\int\Big[\sum_{i=1}^2\r_{x_ix_j}u_{x_kt}\cdot u_{x_ix_jx_k}+\sum_{i=1}^2\r_{x_i}u_{x_kt}\cdot u_{x_ix_j x_jx_k}\\
\di\qquad\qquad\qquad-\sum_{i=1}^2\r_{x_i}u\cdot\nabla u_{x_jx_k}\cdot u_{x_ix_jx_k}-\sum_{i=1}^2\r u_{x_i}\cdot\nabla u_{x_jx_k}\cdot u_{x_ix_jx_k}\Big]dx\\
\end{array}
\end{equation}
and
\begin{equation}\label{u32}
\begin{array}{ll}
\di \int\nabla((\mu+\l(\r)){\rm div} u)_{x_jx_k}\cdot\Delta u_{x_jx_k}dx\\
\di =\int\Big[(\mu+\l(\r))|\nabla({\rm div} u)_{x_jx_k}|^2+\Big(\l(\r)_{x_jx_k}\nabla ({\rm div} u)+\l(\r)_{x_j}\nabla ({\rm div} u)_{x_k}+\l(\r)_{x_k}\nabla ({\rm div} u)_{x_j}\\
\di~+\nabla\l(\r)_{x_jx_k} {\rm div} u+\nabla\l(\r) ({\rm div}
u)_{x_jx_k} +\nabla \l(\r)_{x_j}({\rm div} u)_{x_k}+\nabla
\l(\r)_{x_k}({\rm div} u)_{x_j}\Big)\cdot \nabla({\rm div}
u)_{x_jx_k}\Big] dx.
\end{array}
\end{equation}
Then substituting \eqref{u31} and \eqref{u32} into \eqref{u3}, summing over $j,k=1,2$ and using the Cauchy and Young inequalities and the estimates in Sections 3-4, one has
\begin{equation}\label{ru-3}
\begin{array}{ll}
\di\f{d}{dt}\|\sqrt\r\nabla^3u\|_2^2+2\mu\|\nabla^2\Delta u\|_2^2(t)
 \leq C\Big[(\|u\|^2_{H^3}+1)\|(\nabla^3P(\r),\nabla^3\l(\r))\|^2_2+1\Big].
\end{array}
\end{equation}
Next, applying $\partial_{x_ix_jx_k}$, $i,j,k=1,2$, to
$\eqref{CNS}_1$ gives
\begin{equation}\label{Jan-19-1}
\r_{x_ix_jx_kt}+ \partial_{x_ix_jx_k}({\rm div}(\r u))=0.
\end{equation}
Multiplying \eqref{Jan-19-1} by
$\r_{x_ix_jx_k}$ and summing over $i,j,k=1,2$ and then integrating
with respect to $x$ over $\mathbb{T}^2$, one gets that
\begin{equation}\label{rr1}
\begin{array}{ll}
\di\f{d}{dt}\|\nabla^3\r\|_2^2&\di \leq C\|\nabla^3\r\|_2\Big[\|\nabla\r\|_\i\|\nabla^3 u\|_2+\|\nabla^2 u\|_4\|\nabla^2\r\|_4+\|\nabla u\|_\i
\|\nabla^3\r\|_2+\|\r\|_\i\|\nabla^4u\|_2\Big]\\
&\di \leq C\|\nabla^3\r\|_2\Big[\|\nabla^3 u\|_2+\|\nabla u\|_\i
\|\nabla^3\r\|_2+\|\r\|_\i\|\nabla^2\Delta u\|_2\Big]\\
&\di \leq \a\|\nabla^2\Delta
u\|^2_2+C_\a(\|u\|_{H^3}+1)\|\nabla^3\r\|^2_2,
\end{array}
\end{equation}
where $\alpha>0$ is a constant to be determined.

Similarly, one can obtain
\begin{equation}\label{PP1}
\di\f{d}{dt}\|(\nabla^3P(\r), \nabla^3\l(\r))\|^2_2 \leq \a\|\nabla^2\Delta u\|^2_2+C_\a(\|u\|_{H^3}+1)\|(\nabla^3P(\r),\nabla^3\l(\r))\|^2_2.
\end{equation}
Let $\alpha=\frac{\mu}{3}$. It follows from inequalities \eqref{ru-3},
\eqref{rr1} and \eqref{PP1}, that
\begin{equation}\label{du2}
\begin{array}{ll}
\di\f{d}{dt}\|(\sqrt\r\nabla^3u, \nabla^3\r, \nabla^3P(\r),
\nabla^3\l(\r))\|_2^2(t)+\mu\|\nabla^2\Delta u\|_2^2(t)
 \\
 \di\qquad\qquad\qquad \leq C\Big[(\|u\|^2_{H^3}+1)\|(\nabla^3\r, \nabla^3P(\r),\nabla^3\l(\r))\|^2_2+1\Big].
 \end{array}
\end{equation}
Then integrating \eqref{du2} over $[0,t]$ and using the Gronwall's inequality lead to that
$$
  \begin{array}{ll}
    \di\sup_{t\in[0,T]}\Big[\|\sqrt\r \nabla^3u\|_{2}(t)+\|\nabla^3(\r,P(\r),\l(\r))\|_2\Big]+\int_0^T\|\nabla^2\Delta u\|_2^2 dt\leq C.
   \end{array}$$
So the proof of Lemma \ref{new} is completed. $\hfill\Box$

Now we prove other higher regularities  presented in \eqref{h-regu}
of Theorem \ref{theorem}. It follows easily from \eqref{du2} and \eqref{Jan-16-1}
that for any $t_1, t_2\in [0,T]$,
\begin{equation}\label{P-H3+}
\|\sqrt\r\nabla^3u\|_2^2(t_1)-\|\sqrt\r\nabla^3u\|_2^2(t_2)\to 0,
\end{equation}
as $t_1\to t_2$.

Thanks to  Theorem \ref{theorem2}, one has
\begin{equation}\label{D3+}
\r \in C([0,T];H^2(\mathbb{T}^2))\hookrightarrow
C([0,T]\times\mathbb{T}^2).
\end{equation}

It holds that
\begin{equation}\label{H3-1}
\begin{array}{ll}
\di |\|\r \nabla^3u\|^2_2(t_1)-\|\r \nabla^3u\|^2_2(t_2)|=|\int \r^2|\nabla^3u|^2(t_1,x)dx-\int \r^2|\nabla^3u|^2(t_2,x)dx|\\
\di\leq |\int \r(t_1,x)\big[\r|\nabla^3u|^2(t_1,x)-\r|\nabla^3u|^2(t_2,x)\big]dx|+|\int \r|\nabla^3u|^2(t_2,x)\big[\r(t_1,x)-\r(t_2,x)\big]dx|\\
\di\leq \sup_{[0,T]\times\mathbb{T}^2}\r(t,x)~|\int \big[\r|\nabla^3u|^2(t_1,x)-\r|\nabla^3u|^2(t_2,x)\big]dx|\\
\di \qquad\qquad\qquad\qquad\qquad\qquad\qquad+\sup_{t\in[0,T]}\int \r|\nabla^3u|^2(t,x)dx\sup_{x\in\mathbb{T}^2}|\r(t_1,x)-\r(t_2,x)|\\
\di\leq C\Big[\Big|\int \r|\nabla^3u|^2(t_1,x)dx-\int\r|\nabla^3u|^2(t_2,x)dx\Big|+\sup_{x\in\mathbb{T}^2}|\r(t_1,x)-\r(t_2,x)|\Big]\\
\di \rightarrow 0, ~~{\rm as}~~t_1\rightarrow t_2,
\end{array}
\end{equation}
where one has used  \eqref{P-H3+} and \eqref{D3+}.

Moreover, due to the facts that $\r\nabla^3u\in
L^\infty([0,T];L^2(\mathbb{T}^2)), \r\in C([0,T];H^2(\mathbb{T}^2))$
and $u\in C([0,T];H^2(\mathbb{T}^2))$, it follows that $\r\nabla^3u\in
C([0,T]; H^3-w)$ which means that $\r\nabla^3u$ is weakly continuous
with values in $H^3((\mathbb{T}^2))$. This, together with \eqref{H3-1}, leads
to
\begin{equation}\label{Jan-19-3}
\r\nabla^3u\in C([0,T];L^2(\mathbb{T}^2)).
\end{equation}
In a similar way, one can prove that
 \begin{equation}\label{P-H3}
(\r,P(\r))\in C([0,T];H^3(\mathbb{T}^2))\hookrightarrow C([0,T]; C^1(\mathbb{T}^2)).
\end{equation}
Moreover, since $u\in C([0,T];H^2(\mathbb{T}^2))$ by Theorem
\ref{theorem2} and $\r\in C([0,T]; H^3(\mathbb{T}^2))$ by
\eqref{P-H3}, one can prove that for any $t_1, t_2\in [0,T]$,
\begin{eqnarray}
&&\|\nabla^3\r u(t_1,\cdot)-\nabla^3\r u(t_2,\cdot)\|^2_2\to
0,\label{J-19-1}\\[3mm]
&&\|\nabla\r \nabla^2u(t_1,\cdot)-\nabla\r
\nabla^2u(t_2,\cdot)\|^2_2\to
0,\label{J-19-2}\\[3mm]
&&\|\nabla^2\r \nabla u(t_1,\cdot)-\nabla^2\r \nabla
u(t_2,\cdot)\|^2_2\to 0 \label{J-19-3}
\end{eqnarray}
respectively as $t_1\to t_2$. In fact, to prove \eqref{J-19-3}, one
has
\begin{eqnarray*}
&&\int |\nabla^2\r\nabla u(t_1,x)-\nabla^2\r\nabla u(t_2,x)|^2
dx\\[3mm]
&&\le \int|\nabla^2\r(t_1,x)-\nabla^2\r(t_2,x)|^2 |\nabla
u(t_1,x)|^2 dx+\int |\nabla^2\r(t_2,x)|^2|\nabla u(t_1,x)-\nabla
u(t_2,x)|^2 dx\\[3mm]
&&\le
C\|\nabla^2\r(t_1,\cdot)-\nabla^2\r(t_2,\cdot)\|_2^2+\|\nabla^2\r\|^2_4\|\nabla
u(t_1,\cdot)-\nabla u(t_2,\cdot)\|^2_4\\[3mm]
&&\le
C(\|\nabla^2\r(t_1,\cdot)-\nabla^2\r(t_2,\cdot)\|_2^2+\|\nabla^2
u(t_1,\cdot)-\nabla^2 u(t_2,\cdot)\|^2_2)\to 0,
\end{eqnarray*}
as $t_1\to t_2$. Similarly, \eqref{J-19-1} and \eqref{J-19-2} can be
proved. In view of  \eqref{Jan-19-3} and
\eqref{J-19-1}-\eqref{J-19-3}, we have proved that
\begin{equation}\label{claim}
\r u\in C([0,T]; H^3(\mathbb{T}^2)).
\end{equation}
The proof of Theorem \ref{theorem} is completed.

\section*{Acknowledgments}
Parts of this work were done when Y. Wang was a postdoctoral fellow
at the IMS of Chinese University of Hong Kong during the academic
year 2010-2011; Q.S. Jiu was visiting the IMS of The Chinese University
of Hong Kong, and when Q.S. Jiu and Z.P. Xin were visiting the IMS of
National University of Singapore. The authors would like to thank these
institutions for their supports and hospitality.

\end{document}